
\input amssym.def
\input amssym.tex


\def\item#1{\vskip1.3pt\hang\textindent {\rm #1}}

\def\itemitem#1{\vskip1.3pt\indent\hangindent2\parindent\textindent {\rm #1}}

\tolerance=300
\pretolerance=200
\hfuzz=1pt
\vfuzz=1pt


\hoffset=0.6in
\voffset=0.8in

\hsize=5.8 true in 


\vsize=8.5 true in
\parindent=25pt
\mathsurround=1pt
\parskip=1pt plus .25pt minus .25pt
\normallineskiplimit=.99pt

\countdef\revised=100
\mathchardef\emptyset="001F 
\chardef\ss="19
\def\3{\ss}
\def\anf{$\lower1.2ex\hbox{"}$}
\def\frac#1#2{{#1 \over #2}}
\def\>{>\!\!>}
\def\<{<\!\!<}

\def\into{\hookrightarrow}
\def\ssarr{\hbox to 30pt{\rightarrowfill}}
\def\sarr{\hbox to 40pt{\rightarrowfill}}
\def\arr{\hbox to 60pt{\rightarrowfill}}
\def\larr{\hbox to 60pt{\leftarrowfill}}
\def\Arr{\hbox to 80pt{\rightarrowfill}}

{}

\def\ad{\mathop{\rm ad}\nolimits}

\def\Ad{\mathop{\rm Ad}\nolimits}

\def\cone{\mathop{\rm cone}\nolimits}

\def\det{\mathop{\rm det}\nolimits}

\def\End{\mathop{\rm End}\nolimits}

\def\id{\mathop{\rm id}\nolimits} 
\def\im{\mathop{\rm im}\nolimits}

\def\Ind{\mathop{\rm Ind}\nolimits}

\def\rank{\mathop{\rm rank}\nolimits}



\def\Sl{\mathop{\rm Sl}\nolimits}
\def\SO{\mathop{\rm SO}\nolimits}

\def\Spec{\mathop{\rm Spec}\nolimits}

\def\SU{\mathop{\rm SU}\nolimits}
\def\sup{\mathop{\rm sup}\nolimits}

\def\tr{\mathop{\rm tr}\nolimits}

\def\0{{\bf 0}}
\def\1{{\bf 1}}

\def\a{{\frak a}}

\def\b{{\frak b}}

\def\e{{\frak e}}
\def\f{{\frak f}}
\def\g{{\frak g}}
\def\gl{{\frak {gl}}}
\def\h{{\frak h}}

\def\k{{\frak k}}
\def\l{{\frak l}}
\def\m{{\frak m}}

\def\n{{\frak n}}

\def\p{{\frak p}}
\def\q{{\frak q}}

\def\s{{\frak s}}

\def\sp{{\frak {sp}}}

\def\su{{\frak {su}}}
\def\so{{\frak {so}}}
\def\sL{{\frak {sl}}}
\def\t{{\frak t}}
\def\uu{{\frak u}}

\def\z{{\frak z}}

\def\C{{\Bbb C}}

\def\H{{\Bbb H}}

\def\N{{\Bbb N}}

\def\R{{\Bbb R}}

\def\:{\colon}  
\def\.{{\cdot}}
\def\|{\Vert}
\def\bsk{\bigskip}

\def\giantskip{\vskip2\bigskipamount}
\def\gsk{\giantskip}
\def \la {\langle}
\def\msk{\medskip}
\def \ra {\rangle}
\def \res {\!\mid\!\!}

\def\ssk{\smallskip}

\def\bbr{\bigbreak}
\def\giantbreak{\par \ifdim\lastskip<2\bigskipamount \removelastskip
         \penalty-400 \giantskip\fi}

\def\nin{\noindent}
\def\cen{\centerline}
\def\pagebreak{\vskip 0pt plus 0.0001fil\break}
\def\linebreak{\break}

\def\hat{\widehat}

\def\phi{\varphi}
\def\epsilon{\varepsilon}

\def\nin{\noindent}
\def\oline{\overline}

\def\pder#1,#2,#3 { {\partial #1 \over \partial #2}(#3)}
\def\pde#1,#2 { {\partial #1 \over \partial #2}}


\def\subeq{\subseteq}
\def\supeq{\supseteq}

\def\tilde{\widetilde}

\def\up{{\uparrow}}

\font\eightrm=cmr8


\font\smc=cmcsc10
\font\bfone=cmbx10 scaled\magstep1 
\font\bftwo=cmbx10 scaled\magstep2 

\def\qed{{\unskip\nobreak\hfil\penalty50\hskip .001pt \hbox{}\nobreak\hfil
          \vrule height 1.2ex width 1.1ex depth -.1ex
           \parfillskip=0pt\finalhyphendemerits=0\medbreak}\rm}

\def\qeddis{\eqno{\vrule height 1.2ex width 1.1ex depth -.1ex} $$
                   \medbreak\rm}

\def\Lemma #1. {\bigbreak\vskip-\parskip\noindent{\bf Lemma #1.}\quad\it}

\def\Sublemma #1. {\bigbreak\vskip-\parskip\noindent{\bf Sublemma #1.}\quad\it}

\def\Proposition #1. {\bigbreak\vskip-\parskip\noindent{\bf Proposition #1.}
\quad\it}

\def\Corollary #1. {\bigbreak\vskip-\parskip\nin{\bf Corollary #1.}
\quad\it}

\def\Theorem #1. {\bigbreak\vskip-\parskip\noindent{\bf Theorem #1.}
\quad\it}

\def\Definition #1. {\rm\bigbreak\vskip-\parskip\noindent{\bf Definition #1.}
\quad}

\def\Remark #1. {\rm\bigbreak\vskip-\parskip\noindent{\bf Remark #1.}\quad}

\def\Example #1. {\rm\bigbreak\vskip-\parskip\noindent{\bf Example #1.}\quad}

\def\Problems #1. {\bigbreak\vskip-\parskip\noindent{\bf Problems #1.}\quad}
\def\Problem #1. {\bigbreak\vskip-\parskip\noindent{\bf Problems #1.}\quad}

\def\Conjecture #1. {\bigbreak\vskip-\parskip\noindent{\bf Conjecture #1.}\quad}

\def\Proof#1.{\rm\par\ifdim\lastskip<\bigskipamount\removelastskip\fi\smallskip
            \noindent {\bf Proof.}\quad}

\def\Axiom #1. {\bigbreak\vskip-\parskip\noindent{\bf Axiom #1.}\quad\it}

\def\Satz #1. {\bigbreak\vskip-\parskip\noindent{\bf Satz #1.}\quad\it}

\def\Korollar #1. {\bbr\vskip-\parskip\nin{\bf Korollar #1.} \quad\it}

\def\Bemerkung #1. {\rm\bigbreak\vskip-\parskip\noindent{\bf Bemerkung #1.}
\quad}

\def\Beispiel #1. {\rm\bigbreak\vskip-\parskip\noindent{\bf Beispiel #1.}\quad}
\def\Aufgabe #1. {\rm\bigbreak\vskip-\parskip\noindent{\bf Aufgabe #1.}\quad}

\def\Beweis#1. {\rm\par\ifdim\lastskip<\bigskipamount\removelastskip\fi
           \smallskip\noindent {\bf Beweis.}\quad}

\nopagenumbers

\def\date{\ifcase\month\or January\or February \or March\or April\or May
\or June\or July\or August\or September\or October\or November
\or December\fi\space\number\day, \number\year}

\def\title{Title ??}
\def\author{Author ??}

\def\thanks#1{\footnote*{\eightrm#1}}

\def\rightheadline{\hfil{\eightrm\title}\hfil\tenbf\folio}
\def\leftheadline{\tenbf\folio\hfil{\eightrm\author}\hfil}
\headline={\vbox{\line{\ifodd\pageno\rightheadline\else\leftheadline\fi}}}

\def\firstheadline{}
\def\firstfootline{\cen{\rm\folio}}

\def\seite #1 {\pageno #1
               \headline={\ifnum\pageno=#1 \firstheadline
               \else\ifodd\pageno\rightheadline\else\leftheadline\fi\fi}
               \footline={\ifnum\pageno=#1 \firstfootline\else{}\fi}}

\newdimen\dimenone
 \def\checkleftspace#1#2#3#4{
 \dimenone=\pagetotal
 \advance\dimenone by -\pageshrink   
 \ifdim\dimenone>\pagegoal          
   \else\dimenone=\pagetotal
        \advance\dimenone by \pagestretch
        \ifdim\dimenone<\pagegoal
          \dimenone=\pagetotal
          \advance\dimenone by#1         
          \setbox0=\vbox{#2\parskip=0pt                
                     \hyphenpenalty=10000
                     \rightskip=0pt plus 5em
                     \noindent#3 \vskip#4}    
        \advance\dimenone by\ht0
        \advance\dimenone by 3\baselineskip   
        \ifdim\dimenone>\pagegoal\vfill\eject\fi
          \else\eject\fi\fi}


\def\subheadline #1{\nin\bigbreak\vskip-\lastskip
      \checkleftspace{0.7cm}{\bf}{#1}{\medskipamount}
          \indent\vskip0.7cm\centerline{\bf #1}\medskip}

\def\sectionheadline #1{\bigbreak\vskip-\lastskip
      \checkleftspace{1.1cm}{\bf}{#1}{\bigskipamount}
         \vbox{\vskip1.1cm}\cen{\bfone #1}\bsk}

\def\lsectionheadline #1 #2{\bigbreak\vskip-\lastskip
      \checkleftspace{1.1cm}{\bf}{#1}{\bigskipamount}
         \vbox{\vskip1.1cm}\cen{\bfone #1}\msk \cen{\bfone #2}\bsk}

\def\lchapterheadline #1 #2{\bigbreak\vskip-\lastskip\indent\vskip3cm
                       \cen{\bftwo #1} \msk \cen{\bftwo #2} \gsk}
\def\llsectionheadline #1 #2 #3{\bigbreak\vskip-\lastskip\indent\vskip1.8cm
\cen{\bfone #1} \msk \cen{\bfone #2} \msk \cen{\bfone #3} \nobreak\bsk\nobreak}


\newtoks\literat
\def\[#1 #2\par{\literat={#2\unskip.}%
\hbox{\vtop{\hsize=.15\hsize\nin [#1]\hfill}
\vtop{\hsize=.82\hsize\nin\the\literat}}\par
\vskip.3\baselineskip}

\mathchardef\emptyset="001F 
\def\address{Author: \tt$\backslash$def$\backslash$address$\{$??$\}$}

\def\firstpage{\nin
{\obeylines \parindent 0pt }
\vskip2cm
\centerline {\bfone \title}
\gsk
\centerline{\bf\author}

\vskip1.5cm \rm}

\def\addresstwo{}

\def\dlastpage{\par\vbox{\vskip1cm\nin
\line{
\vtop{\hsize=.5\hsize{\parindent=0pt\baselineskip=10pt\nin\address}}
\quad 
\vtop{\hsize=.42\hsize\nin{\parindent=0pt
\baselineskip=10pt\addresstwo}}
\hfill} }}


\def\firstpage{\nin
{\obeylines \parindent 0pt }
\vskip2cm
\centerline {\bfone \title}
\ssk
\centerline {\bfone \titletwo}
\gsk
\centerline{\bf\author}
\vskip1.5cm \rm}
\def\bs{\backslash} 
\def\addots{\mathinner{\mkern1mu\raise1pt\vbox{\kern7pt\hbox{.}}\mkern2mu
\raise4pt\hbox{.}\mkern2mu\raise7pt\hbox{.}\mkern1mu}}

\pageno=1
\def\up#1{\leavevmode \raise.16ex\hbox{#1}}
 at 8truept
 at 8truept
 at 12truept
\chardef\ss="19
\def\3{\ss}
\def\title{Hardy spaces for non-compactly causal symmetric spaces}
\def\titletwo{and the most continuous spectrum}

\def\author{Simon Gindikin$^*$, Bernhard Kr\"otz${}^{**}$ and 
Gestur \'Olafsson${}^{***}$ }
\footnote{}{${}^*$ Supported in part by the NSF-grant DMS-0070816 and
the MSRI}
\footnote{}{${}^{**}$ Supported in part by the NSF-grant DMS-0097314
and the MSRI}
\footnote{}{${}^{***}$ Supported in part by the NSF-grant DMS-0070607 and
the MSRI}

\def\date{October 29, 2001}
\def\Box #1 { \msk\par\nin 
\centerline{
\vbox{\offinterlineskip
\hrule
\hbox{\vrule\strut\hskip1ex\hfil{\smc#1}\hfill\hskip1ex}
\hrule}\vrule}\msk }

\def\address
{Simon Gindikin

Department of Mathematics

Rutgers University

New Brunswick, NJ 08903

USA

{\tt gindikin@math.rutgers.edu}

\bsk 
\bsk

Gestur \'Olafsson

Louisiana State University

Department of Mathematics

Baton Rouge, LA 70803

USA 

{\tt olafsson@math.lsu.edu}

}

\def\addresstwo
{Bernhard Kr\"otz

The Ohio State University 

Department of Mathematics 

231 West 18th Avenue 

Columbus, OH 43210--1174 

USA

{\tt kroetz@math.ohio-state.edu}

}

\firstpage 

\subheadline{Abstract}

Let $G/H$ be a semisimple symmetric space.  Then
the space $L^2(G/H)$ can be decomposed into
a finite sum of series of  representations induced
from parabolic subgroups of $G$. The most
continuous part of the spectrum of $L^2(G/H)$
is the part induced from the smallest possible
parabolic subgroup. In this paper we
introduce Hardy spaces  canonically related
to this part of the spectrum for a
class of non-compactly causal
symmetric spaces. The Hardy space
is a reproducing Hilbert space
of holomorphic functions  on a bounded symmetric domain of tube type, containing $G/H$ as
a boundary component. A boundary
value map is constructed and we show that
it induces a $G$-isomorphism onto
a multiplicity free subspace of full spectrum 
in the most continuous part $L_{\rm mc}^2(G/H)$
of $L^2(G/H)$. We also relate our Hardy space
to the classical Hardy space on the
bounded symmetric domain .

\sectionheadline{Introduction}

When we transfer from harmonic analysis on Riemannian 
symmetric spaces to non Riemannian semisimple symmetric spaces $G/H$  
one of the most important new phenomena is the fact that
different series of unitary representations appear in the 
decomposition of  $L^2(G/H)$
([BS01ab,D98]).
A main objective of harmonic analysis 
is a {\it geometric} realization of those
{\it series} of representations; an idea that
can be traced back to the article [GG77].
The first step in this program was the realization of the 
holomorphic discrete series for a group $G$ of Hermitian type in 
appropriate Hardy spaces in a curved tube in $G_\C$.
This was accomplished independently in [O82, S86]. 
\par There is a natural generalization of the constructions in [GG77, O82, S86]
for Hermitian groups $G$ to
{\it compactly causal} symmetric spaces $G/H$.
In [\'O\O88] the holomorphic discrete series for $G/H$
was constructed and the appropriate Hardy space was defined and investigated 
in [H\'O\O91]. The explicit form of the Plancherel density
and the corresponding projection operators 
for the holomorphic discrete series were determined by the two last 
named authors [K01,K\'O02, \'O00]. 

\par The above mentioned line of work was on the
{\it discrete} part of the spectrum. In this article we
give for the first time a geometric realization
of the {\it most continuous} part of the spectrum for
a class of non-compactly causal symmetric spaces.

\par  For compactly causal symmetric spaces $G/H$  one has a rich
and
well understood complex geometry: There exist  $G$-invariant tubes
in $G_\C/H_\C$ which have $G/H$ as Shilov boundary; secondly, most
of the compactly causal symmetric spaces can be 
realized as an open dense orbit in the Shilov-boundary of a
bounded Hermitian symmetric space of tube type [Ber96, Bet97, \'O\O99].
This enables us to use complex geometry, holomorphic representations
of semigroups, and the well understood structure and harmonic
analysis of bounded symmetric domains of tube type. 
For non-compactly 
causal symmetric spaces it was believed for a long time 
that  a similar  complex geometrical picture does not exist. 
In [G98] it was 
conjectured that appropriate tubes (related to complex crowns
of Riemannian symmetric spaces) might also exist for non-compactly 
causal symmetric spaces.  It was also conjectured 
that some Hardy spaces on these 
tubes are connected with the  most continuous spectrum.

\par In the last year substantial progress was achieved on the
geometric part of this program. 
This gives us now the possibility 
to start the corresponding  analytical investigations -- 
the subject proper of this paper.

\par The {\it complex crown} $\Xi$ of a 
Riemannian symmetric space $G/K$,  first studied in [AG90], is a certain
open  connected $G$-invariant Stein domain in $G_\C/K_\C$ 
which contains $G/K$ as a totally real submanifold 
(cf.\ Section 1 for the definition). The complex crown 
is universal in many ways. It can be used to
parametrize the compact cycles in complex flag manifolds
(cf.\ [GM01]). Further it has the property 
that the eigenfunctions for the algebra of
$G$-invariant differential operators on $G/K$ extends holomorphically
to $\Xi$ [KS01a]. For a more detailed account on the complex 
crown $\Xi$ we refer to the introduction of [GK02a]. 
\par By the recent work of the two first named authors [GK02b]
it is  known that all non-compactly causal symmetric spaces $G/H$ 
appear in the {\it distinguished boundary} $\partial_d\Xi \subeq G_\C/ K_\C$
of $\Xi$. Furthermore
the non-compactly causal spaces are the only irreducible symmetric spaces
which can be realized as $G$-orbits in $\partial_d\Xi$.

\par In this paper we consider only those complex crowns
that admit a realization as a symmetric space of Hermitian 
type. By this we mean that there 
exists a Hermitian group $S\supeq G$ with maximal 
compact subgroup $U\supeq K$ such that $\Xi$ is $G$-biholomorphic 
to $S/U$. 
For these cases the distinguished  boundary 
$\partial_d\Xi=G(z_1)\simeq G/H$ is exactly one $G$-orbit
and, $\partial_d\Xi$ is equal to the Shilov boundary 
$\partial_s\Xi$ of $\Xi$ in $G_\C/ K_\C$ [GK02b].
Under the additional assumption 
that $G/H$ is irreducible and symmetric, hence 
a non-compactly causal symmetric space, we arrive at the following 
list:
\gsk
\centerline{\bf Table I}
\centerline{\bf  $\Xi\simeq S/U$ and  $\partial_d\Xi\simeq G/H$ is symmetric }

$$\vbox{\tabskip=0pt\offinterlineskip
\def\tablerule{\noalign{\hrule}}
\halign{\strut#&\vrule#\tabskip=1em plus2em&
\hfil#\hfil&\vrule\hfil#\hfil&\hfil#\hfil &\vrule\hfil#\hfil&
\hfil#\hfil&\vrule\hfil#\hfil \tabskip=0pt\cr\tablerule
&&\omit\hidewidth $\g={\rm Lie}(G)$\hidewidth&& 
\omit\hidewidth $\h={\rm Lie}(H)$\hidewidth&& 
\omit\hidewidth $\s={\rm Lie }(S)$\hidewidth& \cr\tablerule
 && $\sp(n,\R)$ && $\gl(n,\R)$ && $\sp(n,\R)\oplus\sp(n,\R)$
&\cr\tablerule
&& $\su(n,n)$  && $\sL(n,\C)\oplus\R$ &&  $\su(n,n)\oplus\su(n,n)$ &
 \cr\tablerule
 && $\so^*(4n)$  && $\sL(n,\H)\oplus\R $&& $\so^*(4n)\oplus\so^*(4n)$ & 
\cr\tablerule
&& $\so(2,n)$ && $\so(1,n-1)\oplus\R$  && $\so(2,n)\oplus\so(2,n)$
&\cr\tablerule
&& $\e_{7(-25)}$ && $\e_{6(-26)}\oplus\R$  &&
$\e_{7(-25)}\oplus\e_{7(-25)}$ 
& \cr\tablerule
&& $\so (1,n)$&&$\so (1,n-1)$&&$\so (2,n)$&\cr\tablerule
 && $\sp(n,n)$ && $\sp(n,\C)$&& $\su(2n,2n)$ &\cr\tablerule
 && $\sp(n,\C)$ && $\sp(n,\R)$ && $\sp(2n,\R)$
&\cr\tablerule}}
$$
  
\par We call the triples $(\s,\g,\h)$ from the above list 
{\it causally symmetric triples.} They can be defined axiomatically 
starting from a Hermitian Lie algebra $\s$ of tube type 
with two commuting involutions $\tau$ and $\sigma$ such that: 
\ssk\nin \ \  (CST1) $\s^\tau=\g$ and $(\s,\tau)$ is compactly causal. 
\par\nin \ \  (CST2) $\g^\sigma=\h$ and  $(\g,\sigma\res_\g)$
is non-compactly causal.

\par Following \'E.\ Cartan the Hermitian symmetric space $S/U$ 
admits a canonical embedding in the dual compact Hermitian symmetric
space and  can be realized as a bounded circled domain in $\C^n$.
We will use Harish-Chandra's 
construction of the realization of $S/U$  as 
a bounded symmetric domain ${\cal D}$. One can
also realize $S/U$ as an affine homogeneous tube domain 
$T_\Omega=\R^n+i\Omega$.  The $\Xi$-realization of $S/U$ we use 
is different from the two mentioned above -- it lies inside 
of $G_\C/K_\C$. 
The different realizations of $S/U$ give
different  Shilov boundaries: A compact Shilov boundary 
$\partial_s{\cal D}$ for ${\cal D}$ which is a single $S$-orbit, and 
a vector space $\R^n$ for the Shilov boundary  of the tubes $T_\Omega$. 
Finally  the Shilov boundary for $\Xi$ is a single $G$-orbit.

\par Hardy spaces are Hilbert spaces of
holomorphic functions which have $L^2$-boundary values on the Shilov 
boundary. Since the three realizations mentioned above  have
different Shilov boundaries, this necessarily leads to 
different definitions of Hardy spaces. 
In the first two realizations the constructions of Hardy spaces 
${\cal H}^2({\cal D})$ and ${\cal H}^2(T_\Omega)$ are well known.
These two Hardy spaces are isomorphic as irreducible 
$S$-modules. We refer to them as {\it classical Hardy spaces}. 
 
\par Our definition of the Hardy space on $\Xi$ is new.
We use a certain minimal semigroup $\Gamma\subeq S_\C$ of elements $\gamma$ with the 
property $\gamma^{-1}{\cal D}\subseteq {\cal D}$ (cf.\ Section 5 for the definition of $\Gamma$). 
Since $\Xi$ is biholomorphic to ${\cal D}$, we obtain 
an action of $\Gamma^{-1}$ on $\Xi$ by compressions. 
We  define the {\it Hardy space} ${\cal H}^2(\Xi)$
on $\Xi$ by 

$${\cal H}^2(\Xi)\:=\{ f\in {\cal O}(\Xi)\: \|f\|^2\:=
\sup_{\gamma\in \Gamma} \int_{G/H} |f(\gamma^{-1} gz_1)|^2\ d\mu_{G/H}(gH)
<\infty\}\ ,$$
where $z_1\in \partial_d\Xi$ is such that $G/H\simeq Gz_1$.
It turns out that ${\cal H}^2(\Xi)$ is a Hilbert space of
holomorphic functions on $\Xi$ with a natural unitary 
action $L$ of $G$ by left translations in the arguments. In particular, 
${\cal H}^2(\Xi)$ admits a reproducing kernel $K_\Xi(z,w)$, the 
so-called {\it Cauchy-Szeg\"o kernel}, holomorphic in the 
first and antiholomorphic in the second variable.

\msk One of our tools is to compare ${\cal H}^2(\Xi)$
with the classical Hardy space  
${\cal H}^2({\cal D})$ which is an 
irreducible $G$-spherical unitary highest weight module 
for the  group $S$. The action of $S$ on ${\cal H}^2({\cal D})$
is given through a cocycle representation 

$$(\pi_h(s)f)(z)=J_h(s^{-1},z)^{-1} f(s^{-1}z)\qquad (f\in {\cal H}^2({\cal D}), 
s\in S, z\in {\cal D}). $$
Theorem 5.7 now  reads:

\msk\nin{\bf Theorem A.} {\it There exists an explicit zero free 
holomorphic function $\psi$ on $\Xi\simeq {\cal D}$ such that 
the mapping 
$$\Psi\: (\pi_h\res_G,{\cal H}^2({\cal D}))\to (L, {\cal H}^2(\Xi)), 
\ \ f\mapsto \psi f$$ 
is a $G$-equivariant isomorphism of Hilbert spaces.}

\msk The fact that $G/K$ is a totally real
submanifold of $S/U$ can now be used to
construct the generalized Segal-Bargmann transform, which
is a unitary $G$-isomorphism $L^2(G/K)\to {\cal H}^2(\Xi )$,
 [D\'OZ01,Ne99,\'O00,\'O\O96,Z01]. This implies that ${\cal H}^2(\Xi )$
is a direct integral of principal series representations parametrized by
$i\a^*/{\cal W}$.
Combining Theorem 5.5 and Theorem 5.8  we come to our 
main  result:

\msk\nin{\bf Theorem B.} {\it The boundary value mapping 
$$b\: {\cal H}^2(\Xi)\to L^2(G/H), \ \ f\mapsto 
\big (gH\mapsto \lim_{\gamma\to \1\atop \gamma\in \Gamma} 
f(\gamma^{-1} gz_1)\big)$$
is a $G$-equivariant isometric embedding. Furthermore
the image of $b$ is a multiplicity free full subspace 
of the most continuous spectrum $L^2(G/H)_{\rm mc}$ 
(cf.\ Definition 4.1).}

\msk Our final result, stated in Theorem 5.9, is the identification
of the reproducing kernel $K_\Xi$ using the fact that
the ${\cal H}^2(\Xi)$ is $G$-isomorphic to ${\cal H}^2({\cal D})$. This
isomorphism clearly carries $K_\Xi$ into the Szeg\"o kernel $K_h$ for
${\cal H}^2({\cal D})$.
The same idea was already used for the compactly causal
case in [B\'O01,\'O\O99].

\msk\nin{\bf Theorem C.} {\it The Cauchy-Szeg\"o kernels $K_\Xi$ 
and $K_h$ are related through
$$K_\Xi(z,w)= \psi(z)\oline {\psi(w)} K_h(z,w). $$}

\msk To illustrate Theorem C consider the example where $S=G\times G$.
Then $\Xi$ is $G$-biholomorphic to ${\cal D}_G\times {\cal D}_G^{\rm opp}$
where ${\cal D}_G$ is the Harish-Chandra realization of $G/K$ and 
${\cal D}_G^{\rm opp}$ refers to ${\cal D}_G$ equipped with the 
opposite complex structure. Denote by $K_h^G(z,w)$ the  
Cauchy-Szeg\"o kernel for the classical Hardy space ${\cal H}^2({\cal D}_G)$.
Then Theorem C says that for all $ (z_1,w_1), (z_2,w_2)\in \Xi\simeq {\cal D}_G\times {\cal D}_G^{\rm opp}$
we have

$$K_\Xi\big((z_1,w_1), (z_2, w_2)\big)={K_h^G(z_1,z_2)K_h^G(w_2,w_1)
\over K_h^G(z_1,w_1) K_h^G(w_2,z_2)}
\ .$$
Note that the formula of $K_\Xi$ is explicit, since 
the classical Hardy space kernels $K_h^G$ are  well known.  

\ssk In the Plancherel decomposition of $L^2(G/H)$ the most 
continuous spectrum $L_{\rm mc}^2(G/H)$ corresponds to one
series of representations, the one which is induced off from a
minimal parabolic subgroup [BS97].  However, the space $L_{\rm mc}^2(G/H)$ is 
not  multiplicity free and we encounter the next problem to 
separate the different multiplicities in $L_{\rm mc}^2(G/H)$
geometrically. Theorem B implies an embedding
$${\cal H}^2(\Xi)\oplus{\cal H}^2(\Xi^{\rm opp})\into 
L^2(G/H)_{\rm mc}$$
which realizes a multiplicity {\it two} subspace of the most
continuous spectrum of $G/H$. Presently it is unclear how to obtain 
the other multiplicities in $L^2(G/H)_{\rm mc}$ missed 
by ${\cal H}^2(\Xi)\oplus{\cal H}^2(\Xi^{\rm opp})$. 
In this context notice the striking similarity with the picture 
for discrete series representations. Only two out of the discrete 
series, the holomorphic and the antiholomorphic, can be realized 
as a Hardy space.

\msk We would like to mention, as we know through several talks, 
that J.\ Faraut is  working on another construction of Hardy spaces 
for certain non-compactly causal symmetric spaces. Also we would 
like to mention that our results in Section 2 are related 
to unpublished work of Y.\ Neretin on explicit branching 
formulas for tensor products of highest weight representations.

\ssk Our paper is organized as follows: 
\msk 
\par 1. Non-compactly causal symmetric spaces and boundary 
components of complex domains
\par 2. Continuous branching of spherical highest weight 
representations
\par 3. Analytical and geometrical constructions on $\Xi$
\par 4. The classical Hardy space inside $L^2(G/H)_{\rm mc}$
\par 5. The Hardy space on $\Xi$
\par Appendix A: parameter calculations
\par Appendix B: Structure theory for causally symmetric triples

\msk It is our pleasure to thank the MSRI, Berkeley, for its hospitality 
during the {\it Integral geometry program} where this 
work was accomplished. We also would like to thank the referee for his careful screening of the 
manuscript and his useful suggestions concerning the readability of the paper. 
 
\pagebreak
\subheadline{}
\centerline{\bfone 1. Non-compactly causal symmetric spaces}
\centerline{\bfone  and  boundary
components of complex domains} 
\gsk

The purpose of this section is mostly of preliminary nature. 
We recall the  basic facts concerning causal 
symmetric spaces. This material is standard and can be found 
in the monograph [H\'O06]. We then turn to the discussion of the complex crown
$\Xi$ of a Riemannian symmetric space $G/K$. We explain some of the new
structural results, in particular those about the distinguished (Shilov) 
boundary $\partial_d \Xi$ of $\Xi$ (cf.\ [GK02a]). Subsequently we
provide the list  of those $\Xi$ which are $G$-isomorphic to a
tube domain $S/U\supeq G/K$
and which have Shilov boundary
$\partial_d\Xi$ isomorphic  to a non-compactly causal symmetric
space $G/H$ (cf.\ Table I in the introduction). This will be  the
class of $\Xi$'s most relevant for this paper. The
section is concluded with a discussion of the
compactification of $\partial_d\Xi$.

\subheadline{Causal symmetric Lie algebras}

\par
Let $\s$ denote  a semisimple real Lie
algebra and $\s_\C$ its complexification. We choose a Cartan involution $\theta$ on $\s$ and write 
$\s=\uu\oplus \p_*$ for the associated Cartan decomposition with 
$\uu$ the maximal compact subalgebra fixed by $\theta$. 

\par In the sequel $\tau\: \s\to\s$ will  denote a non-trivial involution on $\s$ 
commuting with $\theta$. Write $\s=\g\oplus\q_*$ 
for the $\tau$-eigenspace 
decomposition corresponding to the $\tau$-eigenvalues $+1$ and $-1$. 
Then $\g$ is reductive and $\theta\res_\g$ is a
Cartan involution of $\g$. With  $\k =\g\cap\uu$ and 
$\p =\g \cap \p_*$ we obtain the Cartan decomposition $\g=\k\oplus \p$.

\par The symmetric pair
$(\s ,\g)$ is called {\it irreducible} if the only
$\tau$-invariant ideals in $\s$ are the trivial ones,
$\{0\}$ and $\s$. In this case either $\s$ is
simple or $\s \simeq \g\oplus \g$, with $\g$ simple, and $\tau (X,Y)=
(Y,X)$ the flip. 

\par  On the group level we denote by $S$ and $S_\C$
connected Lie groups with Lie algebra $\s$ and $\s_\C$
respectively.
We will assume
-- if not otherwise stated -- that $S\subseteq S_\C$. 
In addition we will require that $\tau$ exponentiates to an involution 
on $S_\C$, again denoted by $\tau$. By $G$ we denote the connected
subgroup of $S$ with Lie algebra $\g$. Finally we call the symmetric space 
$S/G$ {\it irreducible} if $(\s,\g)$ is irreducible.

\Definition 1.1. Let $C$ be an open convex subset of
$\q_*$. Then $C$ is called {\it hyperbolic}, if for all $X\in C$
the operator $\ad(X)$ is semisimple with real eigenvalues. 
We call $C$ {\it elliptic}, if all operators $\ad(X)$, $X\in C$
are semisimple and with imaginary spectrum. 
\qed

We recall some facts on causal symmetric spaces (c.f. [H\'O96],
Chapter 3):

\Definition 1.2. {\bf (Causal symmetric spaces})  The symmetric space $S/G$ is called {\it causal}
if there exists a non-empty  open $G$-invariant convex cone $C$
, containing no affine lines, in $\q_*$.
\qed

There are two different types of causal symmetric spaces,
the non-compactly causal symmetric spaces (NCC) and  the compactly
causal symmetric spaces (CC). In addition, there is the intersection of
those two classes, the Cayley type symmetric spaces (CT).

\Definition 1.3. {\bf (NCC)} Assume that $S/G$ is an irreducible causal symmetric space. Then
the following two conditions are equivalent:

\item{(a)} There exists a non-empty $G$-invariant open hyperbolic cone $C
\subseteq \q_*$ which contains  no affine lines;

\item{(b)} There exists an element $T^0\in \q_*\cap \p_*$, $T^0\not= 0$,
which is fixed by $K$.

If one of those equivalent conditions are satisfied, then $S/G$ is called
{\it non-compactly causal}.
\qed

\Definition 1.4. {\bf (CC)} Assume that $S/G$ is an irreducible causal symmetric space. Then
the following two conditions are equivalent:

\item{(a)} There exists a non-empty  $G$-invariant open elliptic cone $C
\subseteq \q_*$ which contains no affine lines;

\item{(b)} There exists an element $X^0\in \q_*\cap \uu$, $X^0\not= 0$,
which is fixed by $K$.

If one of those equivalent conditions are satisfied, then $S/G$ is
called {\it compactly causal}.
\qed

\Definition 1.5. {\bf (CT)} Assume that $S/G$ is an irreducible causal
symmetric space. Then $S/G$ is called a {\it symmetric space of
Cayley type}, if it is both non-compactly causal and compactly
causal.
\qed

\Remark 1.6. (a) The elements $T^0$ and $X^0$ in Definition 1.3 
and Definition 1.4 are
unique up to multiplication by scalar. If $S/G$ is NCC then we normalize
$T^0$ such that ${\rm ad}(T^0)$ has spectrum $\{0,1,-1\}$. The eigenspace
corresponding to $0$ is exactly $\s^{\theta\tau}=\k \oplus
(\p_*\cap \q_*)$. If $S/G$ is
CC, then we normalize $X^0$ such that the spectrum of
${\rm ad}(iX^0)$ is $\{0,1,-1\}$. In this case the zero eigenspace
is exactly $\uu$.

\par\nin (b) If $S/G$ is compactly causal, then $\uu$ has a non trivial
center $\z (\uu )$ and $\z (\uu )\cap \q_*=\R X^0$. If $S$ is simple then $\z (\uu )
=\R X^0$. If $\s \simeq \g\oplus \g$. Then $X^0 = (X^{00} ,- X^{00} )$
with $X^{00}$ central in $\k $ and
$\z (\uu )=\R X^0\oplus \R (X^{00} ,X^{00} )$. 
\qed

\par
Denote the complex linear extension of $\tau$ to $\s_\C$
by $\tau$.
The $c$-dual $(\s^c,\tau^c)$ of $(\s,\tau)$ is defined 
by $\s^c=\g\oplus i\q_*\subeq \s_\C$ with involution
$\tau^c=\tau|_{\s^c}$. Notice that the $c$-dual of $(\s^c,\tau^c)$
is $(\s ,\tau)$. Let $S^c$ denote the analytic subgroup
of $S_\C$ with Lie algebra $\s^c$. Then $G=  (S\cap S^c)_0$,
where the subscript ${}_0$ denotes the connected
component containing $\1$.
We recall the following fact:

\Proposition 1.7. Assume that $S/G$ is an irreducible symmetric space.
Then the following holds:

\item{(i)} $S/G$ is non-compactly causal if and only if $S^c/G$ 
is compactly causal.

\item{(ii)}  $S/G$ is of Cayley type if and only if
$S/G\simeq S^c/G$.

\item{(iii)} Suppose that $S/G$ is irreducible and causal. Then
$S/G$ is of Cayley type if and only if $G$ is not simple. In that
case $\z (\g )$ is one dimensional and contained in $\p$. 
\qed

\Example 1.8. {\bf (The group case)} Let $G$ be a connected semisimple Lie group,
$S=G\times G$, and $\tau (a,b)=(b,a)$. Then $S/G\simeq G$ is
compactly causal if and only if $G$ is a group of Hermitian type. In this case $S_\C=G_\C\times G_\C$ and
$S^c/G\simeq G_\C/G$ (cf. [H\'O96], Example 1.2.2).
\qed

Let us remark that the causal symmetric pairs $(\s, \g)$ are classified and refer 
to the list [H\'O96, Th.\ 3.2.8].

\subheadline{The domain $\Xi$ -- the complex crown of $G/K$}
According to [GK02a] every  non-compactly causal symmetric
space $G/H$ can be realized as a boundary component of the corresponding 
complex crown $\Xi$ of the Riemannian symmetric space $G/K$. 
In the following two subsections we briefly recall the definition of $\Xi$ and explain some of the 
results of [GK02ab] needed for our purpose. 

\par
Let $\g$ be a reductive Lie algebra with complexification $\g_\C$. 
We let $G_\C$ be a connected complex Lie group with Lie algebra $\g_\C$
and let $G$ be the analytic subgroup of $G_\C$ corresponding to $\g$. 
Let $K\subseteq G$ be a maximal compact subgroup and let $\k$ be its  Lie algebra. 
Let $\g=\k\oplus\p$ be the Cartan decomposition corresponding to $\k$.
Finally we choose a maximal abelian subspace $\a\subeq 
\p$. Let $\m =\z_{\k}(\a)$ and
for $\alpha\in \a^*$ define 
$$\g^\alpha =\{X\in \g\:(\forall H\in \a ) [H,X]=\alpha (H)X\}\, .$$
Let $\Sigma =\Sigma (\g ,\a )$ be the corresponding
set of (restricted) roots, i.e., $\Sigma =\{\alpha \in \a^*\setminus \{0\}\:
\g^\alpha \not= \{0\}\}$.
Then
$$\g=\a\oplus \m \oplus\bigoplus_{\alpha\in \Sigma} \g^\alpha\, .$$ 
We write $A$ respectively $N$ for the analytic subgroup of $G$ corresponding to 
$\a$ respectively $\n$.  
We also need the {\it Weyl
group} of $\Sigma (\a ,\g )$,  ${\cal W} \:=N_K(\a)/ Z_K(\a)$.

Let
$$\Omega=\{ X\in\a\:(\forall \alpha\in \Sigma)\  |\alpha(X)|<{\pi\over
2}\}\ .$$
Then the {\it complex crown} of the Riemannian symmetric space 
$G/K$ is defined by 
$$\Xi\:= G\exp(i\Omega) K_\C/ K_\C\, .$$  
Note that $\Xi$ is an open domain in $G_\C/ K_\C$ and that 
$G$ acts properly on it (cf.\ [AG90]). 
Write $\partial \Xi$ for the topological boundary of $\Xi$
and notice that $G\exp(i\partial \Omega)K_\C/ K_\C\subeq \partial
\Xi$ (cf.\ [AG90]).

\subheadline{The distinguished boundary of $\Xi$}

Write $\oline\Xi$ for the closure of $\Xi$ in $G_\C/ K_\C$, and
note that $\oline \Xi =\Xi\amalg\partial \Xi$, where $\amalg$ denotes
disjoint union.  
Write $\partial_e \Omega$ for the extreme points of 
the compact convex set $\oline{\Omega}$. Recall from 
[GK02a] that $\partial_e {\Omega}$ is a finite union 
of ${\cal W}$-orbits:
$$\partial_e\Omega={\cal W}(X_1)\amalg\ldots\amalg{\cal W}(X_n)\ .$$ 
We define the {\it distinguished boundary  of $\Xi$ in $G_\C/K_\C$}
by 
$$\partial_d\Xi\:=G\exp(i\partial_e\Omega)K_\C/K_\C\ .$$ 
It was shown in [GK02a, Th.\ 2.3] that $\partial_d\Xi$
has the property that
$$\sup_{z\in \Xi} |f(z)|=\sup_{z\in \partial_d\Xi} 
|f(z)|\leqno(1.1)$$
for all bounded holomorphic functions $f$ on $\Xi$ which 
continuously extend to $\oline \Xi$. Moreover,  if $\partial_d\Xi$
is connected, i.e., if $\partial_e\Omega={\cal W}(X_1)$ is 
a single orbit, then $\partial_d\Xi$ is minimal with respect to the 
property (1.1). In this case  we call $\partial_d\Xi$ the {\it Shilov boundary}
of $\Xi$ and denote it by $\partial_s\Xi$.

\par Write $z_j\:=\exp(iX_j)K_\C \in \partial_d\Xi$ and write 
$H_j$  for the isotropy subgroup of $G$ in $z_j$. Then we have 
the following (cf. [GK02a, Th.\ 3.6, Th.\ 3.26]):

\Theorem 1.9. Assume that $G_\C$ is simply connected. Then, with 
the notation introduced above, the 
following assertions hold:
\item{(i)} $\partial_d\Xi\simeq G/H_1\amalg\ldots\amalg G/H_n$.
\item{(ii)} If one of the boundary components 
$G/H_j$ is an irreducible  symmetric space, then it is a
non-compactly causal symmetric space. Moreover every 
non-compactly causal symmetric space $G/H$ attached 
to $G$ appears $G$-locally in the distinguished boundary 
of $\partial_d\Xi$ of $\Xi$. \qed

\subheadline{Realization of $\Xi$ as tube domains}

We consider now the case where $\Xi$ has a realization as a
tube domain. The material in this subsection is taken from [KS01b]. 
We go back to our notations from the first subsection,
as we will be considering causally symmetric triples $\h\subseteq \g\subseteq \s$.
Thus throughout this section $(\s,\tau)$ denotes a compactly
causal symmetric Lie algebra. We require $S_\C$ to be simply connected
and write $S, G, U$ and $K$ for the analytic subgroups of $S_\C$ 
with Lie algebras $\s, \g, \uu$ and $\k$. 
Note that the embedding 
$$G/K\into S/U, \ \ gK\mapsto gU$$
realizes $G/K$ as a totally real submanifold of the 
Hermitian symmetric space $S/U$.

\par Recall that $\z(\uu)\cap \q_*=\R X^0$ is one dimensional
and that we can normalize $X^0$ such
that $\Spec(\ad iX^0)=\{0,1,-1\}$.
We let $\t\subeq \uu$ be a compact $\tau$-stable 
Cartan algebra of $\s$ and note that 
$\z(\uu)\subeq \t$. Write $\Delta=\Delta(\s_\C,\t_\C)$ for 
the root system with respect to $\t_\C$. 
Call a root $\alpha $ {\it compact} if $\alpha\res_{\z(\uu)}=0$
and {\it non-compact} otherwise. Thus $\alpha $ is non-compact
if and only if $\g^\alpha_\C\subset \p_\C$. 
Denote by $\Delta_c$ the set of
compact roots and by $\Delta_n$ the set of non-compact roots. 
Let $\Delta_n^+:=\{\alpha \in \Delta\: \alpha (iX^0)=1\}$
and fix a positive system $\Delta^+$ such that
$\Delta_n^+\subseteq \Delta^+$.
Write $\p^{\pm}=\oplus_{\alpha\in \Delta_n^+}\g^{\pm \alpha}_\C\subseteq \p_\C$.
Then we have the triangular decomposition of $\s_\C$:
$$\s_\C=\p^+\oplus\uu_\C\oplus\p^-\ .$$
Let $P^\pm =\exp (\p^\pm)$. Let $U_\C\subseteq S_\C$ be
the complexification of $U$.
Then it is well known that  $S\subeq P^+ U_\C P^-$
and we have the Borel 
embedding 
$$S/U\into S_\C/ U_\C P^-\ .$$
The map $P^+\times U_\C \times P^-\ni (p^+,u,p^-)\mapsto p^+up^-
\in P^-U_\C P^+$
is a diffeomorphism onto the open dense subset $P^+ U_\C P^-\subseteq S_\C$.
We write $s=p^+(s)\kappa (s)p^-(s)$ for the unique decomposition
of an element $s\in P^+U_\C P^-$. Then
we have the Harish Chandra realization of 
$S/U$ as a bounded symmetric domain ${\cal D}\subeq \p^+$
given by
$$S/U\ni sU\mapsto \log (p^+(s))\in \p^+\ , $$
where $\log =(\exp|_{\p^+})^{-1}:P^+\to \p^+$.

\par Write $\Sigma=\Sigma(\g,\a)$ for the restricted root system 
with respect to $\a$ and $\hat\Sigma=\Sigma(\s,\a)$ for 
the double restricted root system. If $\Sigma$ is an abstract 
root system on $\a$, then  we also write $\Omega=\Omega(\Sigma)$ to
indicate 
the dependency of $\Omega$ on $\Sigma$. Define 

$$\Xi_0\:=G\exp(i\Omega(\hat\Sigma))K_\C/ K_\C\ ,$$
and notice that $\Xi_0\subeq \Xi$ since $\Sigma\subeq \hat\Sigma$.

\Theorem 1.10. (cf.\ {\rm [KS01b, Sect.\ 2]}) The embedding $G/K\into S/U$ 
extends to a $G$-equivariant 
biholomorphism
$$\Phi\: \Xi_0\to S/U\subeq S_\C/ U_\C P^-\ .$$
Furthermore the following statements are equivalent:
\item{(1)} $\Xi=\Xi_0$. 
\item{(2)} $\Omega(\hat \Sigma)=\Omega(\Sigma)$. 
\item{(3)} $\rank_\R (\g)={1\over 2}\rank_\R(\s)$. 
\item{(4)} $\Sigma$ is of type $C_n$ or $BC_n$
for $n\geq 2$ or $(\g ,\s )$ is one of $(\su (1,1),\su (1,1)\oplus \su (1,1))$
or $(\so (1,n), \so(2,n))$. \qed

The non-compactly causal symmetric spaces $G/H$ 
which are singled out by the additional requirement that 
$\Xi=\Xi_0\simeq S/U$ are the ones that are  most important for 
the rest of this article. We will
therefore discuss them in more details now. 
The non-compactly causal symmetric
spaces that satisfy the additional requirement $\Xi=\Xi_0$ 
are exactly those where $\Sigma$ is of type $C_r$, 
say 
$$\Sigma=\{ {1\over 2}(\pm\gamma_i\pm\gamma_j)\: 1\leq i,j\leq r\}\bs\{0\}\ .$$
Define $Y^j\in \a$ by $\gamma_i(Y^j)=2\delta_{ij}$
and $Y^0=\frac{1}{2}(Y^1+\ldots +Y^r)$. In this case 
$n=1$, $X_1=\frac{\pi}{2}Y^0$ and
$$\partial_e\Omega={\cal W} ({\pi\over 2} Y^0)\, .$$
Hence with $z_1\:=\exp(i{\pi\over 2}Y^0)K_\C$ 
we have that $\partial_s\Xi=G(z_1)$ 
is a connected  $G$-space. Write $H$ for the isotropy 
subgroup of $G$ in $z_1$ and $\h$ for the Lie algebra
of $H$. In cases where ${\rm Ad}(z_1)$ plays the role of a partial Cayley transform
we write ${\bf c}$ for $z_1$. 
It follows from [GK02b]
that 
$$\sigma\: G\to G, \ \ g\mapsto (\Ad({\bf c}^2)\circ\theta)(g)$$
is an involution on $G$ such that
$H=G^\sigma$. Comparing Theorem 1.10 with the list of all causal pairs 
in [H\'O96, Th.3.2.8] we finally arrive at Table I from the introduction. 

\par Note that the first five  pairs $(\g,\h)$ in Table I  
are of Cayley type. 
We call a triple $(\s,\g,\h)$ from Table I  a {\it causally
symmetric triple}. Note 
that we always have 
$$\rank_\R(\h)=\rank_\R(\g)={1\over 2}\rank_\R (\s)\ .$$
We write $\g=\h\oplus\q$ for the $\sigma$-eigenspace decomposition
of $\g$. Note that $\sigma$ extends to an involution 
on $\s$ since $\Spec \ad_\s(Y^0)=\{ -1, 0, 1\}$. Further observe that $\sigma$
commutes with $\tau$. Our choice of the maximal abelian subspace $\a\subeq \p$
is then such that $\a\subeq \q\cap\p$. 

\par Causally symmetric triples can also be defined axiomatically: 

\Definition 1.11. {\bf (Causally symmetric triples)} A triple 
of Lie algebras $(\s,\g,\h)$  with $\h\subeq \g\subeq \s$ is 
called a {\it causally symmetric triple} if the following 
axioms are satisfied: 
\ssk\nin \ \ (CST1) $\s$ is Hermitian and of tube type.
\ssk\nin \ \ (CST2) There exists two commuting involutions $\tau$ and $\sigma$ on 
$\s$ such that:
\itemitem{(1)}  $(\s,\tau)$ is compactly causal and $\g=\s^\tau$. 
\itemitem{(2)}  $(\g,\sigma\res_\g)$ is 
non-compactly causal and $\h=\g^\sigma$. \qed

\ssk If we only require that $G/K$ embeds into $S/U$ 
and that $\partial_d\Xi$ has at least one boundary $G/H_j$ 
component (cf.\ Theorem 1.9) which is a non-compactly causal symmetric 
space, then we arrive at the following more general situation 
(cf.\ [KS01b] and [GK02b]). 

\Theorem 1.12. Suppose that $G/K$ is a totally real
submanifold of $S/U$ and that there exists a symmetric
subgroup $H$ of $G$ such that $G/H$ is a non-compactly
causal symmetric space and is a $G$-orbit in $\partial_d\Xi$.
With $\partial_d\Xi_0=G\exp(i\partial_e \Omega(\hat\Sigma))K_\C /K_\C$
the distinguished boundary of $\Xi_0$ we arrive at 
the following four possibilities:
\item{(i)} $(\s ,\g ,\h)$ is a causally symmetric
triple.
\item{(ii)} $G$ is the
structure group of an Euclidian Jordan algebra, i.e., $\g$ is 
isomorphic to one of the following
$$\gl(n,\R)\qquad \sL(n,\C)\oplus\R\qquad \sL(n,\H)\oplus\R\qquad 
\so(1,n)\oplus\R\qquad \f_{4(-20)}\oplus\R\ ,$$
and $\partial_d\Xi_0=\partial_d\Xi\amalg G/K\amalg G/K$. 
\item{(iii)} $G=\SO (n,n)$,
$\partial_d\Xi_0=G/\SO (n,\C)$,
and
$$\partial_d\Xi = \partial_d \Xi_0\amalg G/
\big(\SO (1,n-1)\times \SO (1,n-1)\big)\, .$$
\item{(iv)} $G=\SO (2n,\C )$,
$\partial_d\Xi_0=G/\SO^*(2n)$, 
and
$\partial_d\Xi =\partial_d\Xi_0\amalg G/\SO (2,2n-2)$.
\qed

\subheadline{The compactification of the Shilov boundary
$\partial_s\Xi$}

For the remainder of this section we assume that $(\s,\g,\h)$ is 
a causally  symmetric triple.

\Proposition 1.13. Assume that $(\s,\g,\h)$ is a causally symmetric 
triple. 
Then the mapping 
$$\iota\: G_\C/ K_\C \to S_\C/ U_\C P^-, \ \ gK_\C\mapsto gU_\C P^-$$
is an open $G_\C$-equivariant embedding. 

\Proof. Since $\s_\C = \g_\C +\frak{u}_\C + \p^-$ it 
is clear that the image of $\iota$ is open. It remains to show that 
$\iota$ is injective, which means that $G_\C\cap U_\C P^-=K_\C$. 
Also denote by $\tau$ and $\theta$ the holomorphic extension of $\tau$
and $\theta$ to an involution on $S_\C$. Note that $S_\C^\tau=G_\C$ 
and $S_\C^\theta=U_\C$ since $S_\C$ was assumed to be simply connected.

\par We have $\tau(P^-)=P^+$. Thus we obtain 
$$G_\C\cap U_\C P^-= G_\C \cap U_\C P^-\cap \tau (U_\C P^-)
=G_\C \cap \big (U_\C P^-\cap U_\C P^+\big)\ .$$
Now the well known fact $U_\C P^-\cap U_\C P^+= U_\C$
implies  that 
$$G_\C \cap U_\C P^-=G_\C \cap U_\C\ .$$
It remains to show that $G_\C\cap U_\C =K_\C$. Now 
we are going to use the fact that the pair $(\s,\g,\h)$ is 
from Table I. Quick inspection shows that $G_\C$ is 
simply connected. Hence $G_\C^\theta=K_\C$ is connected
and so $G_\C\cap U_\C=G_\C\cap S_\C^\theta=K_\C$, concluding 
the proof of the proposition.\qed 

Write $\oline {\cal D}$ for the closure in $\p^+$ of the bounded symmetric
domain ${\cal D}$. Recall that this is also the closure of $S/U$
in $S_\C/ U_\C P^-$. Hence Proposition 1.13 implies that 
the isomorphism $\Phi\: \Xi\to {\cal D}$ from Theorem 1.10 extends 
to a $G$-equivariant embedding 
$$\oline\Phi\: \oline\Xi\to \oline {\cal D}\ .$$
Thus $\oline {\cal D}$ is the natural $G$-invariant 
compactification of $\Xi$. Let
$z_1\:=\exp(i{\pi\over 2}Y^0)K_\C \in \partial_d\Xi$ as before.
It is known that the stabilizer of $z_1$ in $S$ is 
a maximal parabolic subgroup $P$ of $S$. Note that 
$P$ is given by 
$$P={\bf c}U_\C P^-{\bf c}^{-1}\cap S\ .$$
Furthermore the Shilov boundary 
$\partial_s{\cal D}$ of ${\cal D}$ is the $S$-orbit through $z_1$
and 
$$\partial_s{\cal D}= S(x_1)=U(x_1)\simeq U/K\ .$$
{}From the $G$-equivariance of $\oline\Phi$ it hence follows 
that $\oline \Phi(\partial_s\Xi)\subeq \partial_s{\cal D}$.

\Theorem 1.14. The image of the Shilov boundary $\partial_s\Xi$
under the embedding $\oline\Phi$ is open and dense 
in $\partial_s{\cal D}$. In particular, $\partial_s{\cal D}$
is a compactification of $\partial_s\Xi$. 

\Proof. First it is clear that 
$\oline\Phi(\partial_s\Xi)$  is contained in $\partial_s{\cal D}$. Further 
note that $G/H\into S/P$ is an open 
$G$-invariant embedding. Hence $\oline\Phi(\partial_s\Xi)$ is open. 
\par Let us now show that $\oline \Phi(\partial_s\Xi)$ is dense in $\partial_s{\cal D}$. 
For that we realize $\oline \Xi$ in $\oline{\cal D}$
via the embedding $\oline\Phi$. 
If $Y\:=\oline\Phi(\partial_s\Xi)$ is not dense,  then we find a 
continuous function $f$ on $\oline{\cal D}$, holomorphic on ${\cal
D}=\Xi$
which does not attain its maximum on $\oline Y$. 
But this contradicts the fact that $\Xi={\cal D}$ 
and that $\partial_s\Xi$ is the Shilov boundary of $\Xi$.\qed

\sectionheadline{2. Continuous branching of spherical highest weight
representations}

Throughout this section $(\s,\tau)$ denotes a compactly 
causal symmetric Lie algebra with $\tau$-eigenspace 
decomposition $\s=\g\oplus\q_*$. Denote by $S$  a simply
connected Lie group with Lie algebra $\s$ and write
$G\:=S^\tau$ for the fixed point group. Observe that
$G$ is connected since $S$ is assumed  to be simply connected.
\par In this section we discuss the
{\it generalized Segal-Bargmann transform}
introduced in [\'O{\O}85] (see also [D\'OZ01,Ne99,\'O00,Z01]). 
The generalized Segal-Bargmann transform is a unitary $G$-isomorphism
$U_\lambda:(L,L^2(G/K)) \to (\pi_\lambda, {\cal H}_\lambda)$ where ${\cal H}_\lambda$ is
a unitary $G$-spherical highest weight representation of $S$.
In this section we will give a sharp criterion in terms of $\lambda$ for the
existence of $U_\lambda$. In particular this result will 
imply the continuous branching of $\pi_\lambda|_{G}$.

\msk  
Write $U$ for the analytic subgroup of $S$ associated 
to $\uu$ and $U_\C$ for its complexification. We denote by 
$P^\pm$ the analytic subgroup of $S_\C$ with Lie algebras
$\p^\pm$.
Assume for a moment that $S$ is not simply connected and that 
$S\subeq S_\C$. By the Harish Chandra decomposition of $S$ we have 
$S\subeq P^+ U_{\C} P^-$. For $s\in S$ write
$s=p^+(s)\kappa (s)p^-(s)$ with the obvious notation. Recall
that the map $S\ni s\mapsto \log p^+(s)\in \p^+$ gives
the realization of $S/U$ as a bounded symmetric
domain in $\p^+$. Further we have the middle projection
$\kappa :S \rightarrow U_\C$.  By the usual lifting 
argument we obtain a map $\kappa \: S\to U_\C$  for 
a simply connected group $S$. This lifting is unique 
if we  require -- as we will -- that $\kappa(\1)=\1$. We denote by
$X\mapsto \oline X$ the conjugation in $\s_\C$ with
respect to the real form $\s$. Then
$\overline{\p^+}=\p^-$. 

\subheadline{Realization of unitary highest weight representations.}

Our source of reference for the facts collected below is 
[N99, Ch.\ XII]. Let $(\pi_\lambda, {\cal H}_\lambda)$ be a unitary 
highest weight representation of $S$ with highest weight 
$\lambda\in i\t^*$ and with respect to the positive 
system $\Delta^+$. Write $F(\lambda)$ for the finite dimensional 
lowest $U$-type in ${\cal H}_\lambda$, i.e., the space 
which is generated by applying $U$  to a highest weight vector 
$v_\lambda$ of ${\cal H}_\lambda$. 
We now briefly recall the realization of ${\cal H}_\lambda$
in ${\cal O}({\cal D}, F(\lambda))$. For that we write 
$\sigma_\lambda$ for the representation of $U_\C$ on $F(\lambda)$.
Define
$${\cal K}\: {\cal D}\times{\cal D}\to U_\C , 
\ \ (z,w)\mapsto \kappa(\exp(-\oline w)\exp(z))^{-1}$$
$$J\: S\times {\cal D}\to U_\C, \ \ 
(s,z)\mapsto \kappa(g\exp(z))\ .$$
Set
$$K_\lambda=\sigma_\lambda\circ {\cal K}\, \quad {\rm and}\quad 
J_\lambda=\sigma_\lambda\circ J\, .\leqno(2.1)$$ 
Note that 
$K_\lambda$ is holomorphic in the first variable,  
antiholomorphic in the second and satisfies the 
cocyle identity
$$K_\lambda(s(z), s(w))=J_\lambda(s,z) K_\lambda(z,w)
J_\lambda(s,w)^*\leqno(2.2)$$
for all $s\in S$, $z,w\in {\cal D}$. Here $A^*$ denotes the adjoint 
of an operator $A\in B(F(\lambda))=\End(F(\lambda))$. 
\par The relation to $(\pi_\lambda, {\cal H}_\lambda)$ is 
as follows: $\pi_\lambda$ is unitary if and only if 
$K_\lambda$ is positive definite. If this is the case, then 
$(\pi_\lambda, {\cal H}_\lambda)$ can be realized 
in ${\cal O}({\cal D}, F(\lambda))$ as the reproducing 
kernel Hilbert space corresponding to $K_\lambda$. 
So henceforth we will assume that ${\cal H}_\lambda\subeq {\cal O}
({\cal D}, F(\lambda))$. The representation $\pi_\lambda$ is 
then given by 

$$(\pi_\lambda(s)f)(z)=J_\lambda(s^{-1},z)^{-1}f(s^{-1}z)\leqno(2.3)$$
for $f\in {\cal H}$, $s\in S$ and $z\in {\cal D}$.

\subheadline{$G$-spherical unitary highest weight representation}

A unitary highest weight representation $(\pi_\lambda, 
{\cal H}_\lambda)$ is called {\it $G$-spherical} if there
exists a non-trivial $G$-fixed  element $\nu\in {\cal H}_\lambda^{-\infty}$,
where ${\cal H}_\lambda^{-\infty}$
denotes the  module of distribution vectors of 
$(\pi_\lambda, {\cal H}_\lambda)$. 
If $(\pi_\lambda, {\cal H}_\lambda)$ is a unitary 
highest weight representation of $S$, then we call 
$\lambda$ {\it regular} if ${\cal H}_\lambda$ is dense 
in ${\cal O}({\cal D}, F(\lambda))$, or, equivalently 
if $\lambda\res_{\z(\uu)}$ lies on the continuous
halfline of the Wallach set (this is made precise in Appendix A).
We say that $(\pi_\lambda,{\cal H}_\lambda)$ is {\it scalar} if $F(\lambda )=\C$.
This is equivalent to $\lambda\in i\z(\uu)^*$.

\par We now recall some facts concerning $G$-spherical 
unitary highest weight representation. For more
information and an almost complete classification see [KN02].

\msk\nin  \item{$\bullet$} If $(\pi_\lambda, {\cal H}_\lambda)$ is 
$G$-spherical, then $\nu_0\:=\nu\res_{F(\lambda)}$ is
non-zero and $K$-fixed. In particular, $F(\lambda)$ is a $K$-spherical 
module for $U$. 
\item{$\bullet$} Suppose that $ \lambda$ is regular, 
then $(\pi_\lambda, {\cal H}_\lambda)$ is $G$-spherical if and 
only if $F(\lambda)$ is $K$-spherical. 
\item{$\bullet$} Suppose that $ \lambda\in i\z(\uu)^*$ is regular, 
then $(\pi_\lambda, {\cal H}_\lambda)$ is $G$-spherical if and 
only if $\lambda\in i(\z(\uu)\cap\q_*)^*$. In particular, if 
$S$ is simple, then all regular $\lambda\in i\z(\uu)^*$ 
correspond to $G$-spherical representations. 
\item{$\bullet$} If we are in the group case $S=G\times G$, 
then the $G$-spherical unitary highest weight representations
of $S$ are exactly the representations $\pi_\lambda\otimes\pi_\lambda^*$
with $\pi_\lambda$ a unitary highest weight representation of $G$
and $\pi_\lambda^*$ its dual representation. 
\msk 
\subheadline{The generalized Segal-Bargmann transform}

In this subsection $(\pi_\lambda, {\cal H}_\lambda)$ denotes 
a $G$-spherical unitary highest weight representation of $S$.  
Note that this means in particular that there exists a non-trivial
$K$-fixed vector $\nu_0\in F(\lambda)^*$. Write $\nu_0=\la \cdot , v_0 \ra$
where $0\not= v_0\in F(\lambda)$ is a $K$-fixed vector. 
Define 
$$D_\lambda\: G\to B(F(\lambda)), \ \  g\mapsto J_\lambda (g,0)^{-1}=
\sigma_\lambda (\kappa (g))^{-1}$$
and
$$\|D_\lambda\|\: G/K\to \R^+, \ \ gK\to \|D_\lambda(g)\| $$
where $\|\cdot\|$ denotes the operator norm on $B(F(\lambda))$. 
Write ${\cal A}(G/K)$ for the space of analytic functions on 
$G/K$. Define a restriction map
$$R_\lambda :{\cal H}_\lambda \to {\cal A}(G/K), \ \ 
R_\lambda (f)(gK):=\la D_\lambda (g)f(gK), v_0\ra\ .$$
Note that this map is well defined, since $v_0$ is $K$-fixed. 
We write $L$ for the left regular representation of $G$ on 
functions on $G/K$, i.e., $(L(g)f)(xK)=f(g^{-1}xK)$ for $g,x\in G$ and 
$f$ a function on $G/K$. 
We then have the following lemma, see [\'O00], Lemma 3.4:

\Lemma 2.1. Let $(\pi_\lambda, {\cal H}_\lambda)$ be a $G$-spherical 
unitary highest weight representation of $S$. 
\item{(i)} The map $R_\lambda\: {\cal H}_\lambda\to {\cal A}(G/K)$
intertwines $\pi_\lambda\res_G$ with $L$. Moreover, if 
$\lambda\in i\z(\uu)^*$, then $R_\lambda$ is injective. 
\item{(ii)} If $\lambda$ is regular and $\|D_\lambda\|\: G/K\to \R^+$
is square integrable, then $\im R_\lambda\subeq L^2(G/K)$ is dense 
and $R_\lambda\: {\cal H}_\lambda\to L^2(G/K)$ is continuous.\qed 

Assume now that $\|D_\lambda\| \in L^2(G/K)$ and  that $\lambda$ is 
regular. Then $R_\lambda\: {\cal H}_\lambda\to L^2(G/K)$ is 
continuous with dense image. Thus we can consider the polar 
decomposition of the continuous operator $R_\lambda$ given by 
$R_\lambda =U_\lambda P_\lambda$ with 
$P_\lambda =(R_\lambda^* R_\lambda)^{1\over 2}$ and $U_\lambda$ 
a partial isometry. If $\lambda\in i\z(\uu)^*$, then 
$U_\lambda$ is an isometry and 
the unitary $G$-equivariant isomorphism 
$$U^*_\lambda :L^2(G/K)\to {\cal H}_\lambda$$
is called the {\it generalized Segal-Bargmann transform} (cf.\
[\'O00]). 

\par We will now specify the set of parameters for which 
$\|D_\lambda\|$ is square integrable. We will distinguish 
the two cases of $S$ simple and $S=G\times G$. The group 
case is much simpler, since in that case the square integrability 
of $D_\lambda$ can be easily  reduced to known integrals.
We therefore treat the 
group case first since the general proof requires some knowledge 
of advanced structure theory the general reader might not be 
so familiar with.

\subheadline{The square integrability of $\|D_\lambda\|$: group case}

In this section we consider the case where $S= G\times G$. 
We let $(\pi_\lambda, {\cal H}_\lambda)$ be a unitary 
representation of $G$. 
Then the $S$-representation we consider is $\hat\pi_\lambda =
\pi_\lambda\otimes\pi_\lambda^*$. Recall that these are 
precisely the $G$-spherical unitary highest weight representations
of $S$. As explained before we can realize ${\cal H}_\lambda
$ in ${\cal O}(G/K, F(\lambda))$. As before we write $K_\lambda$ for its 
reproducing kernel and $J_\lambda$ for the corresponding 
cocycle.
We recall now some well known facts on representations of 
the holomorphic discrete series on $G$. Write $Z$ for the center of $G$. We say 
$(\pi_\lambda,{\cal H}_\lambda)$ belongs to the {\it relative 
discrete series} if all the matrix coefficients $g\mapsto (\pi_\lambda (g)u,v)$,
$u,v\in {\cal H}_\lambda$, are square integrable modulo 
$Z$. 
By abuse of notation we will for the moment denote by 
$\t$ a compact Cartan algebra of $\g$ and 
$\Delta$ for the corresponding root system. 
We write $\beta$ for the highest root in $\Delta^+$
and  set $\rho\:={1\over 2}\sum_{\alpha\in \Delta^+}\alpha$.

\Lemma 2.2. Let $(\pi_\lambda, {\cal H}_\lambda)$ be a unitary 
highest weight representation of the simply connected Hermitian 
Lie group $G$. Then the following assertions are equivalent:
\item{(1)} $(\pi_\lambda, {\cal H}_\lambda)$ belongs 
to the relative discrete series. 
\item{(2)} The highest weight $\lambda\in i\t^*$ satisfies 
the Harish Chandra condition 
$$\la \lambda+\rho,\beta\ra <0\ .$$
\item{(3)} The integral 
$$\int_{G/K} \|K_\lambda(z,z)\|^{-1}\ d\mu_{G/K}(z)<\infty$$
is finite. Here $\|\cdot\|$ denotes the operator norm on 
$B(F(\lambda))$ and $\mu_{G/K}$ the invariant 
measure on $G/K$. 
\item{(4)} The integral 
$$\int_{G/Z} \|J_\lambda(g,0)^{-1}(J_\lambda(g,0)^{-1})^*\| 
\ d\mu_{G/Z}(gZ)$$
is finite. 

\Proof. (1)$\iff$(2) is the criterion of Harish Chandra (cf.\ [N99, 
Th.\ XII.5.12]). For the
equivalence 
of (2) and (3) see [N99, Th.\ XII.5.6]. Finally, the equivalence of (3)
and (4) is immediate from (2.1) and (2.2).\qed 

As before we realize 
$(\hat\pi_\lambda, {\cal H}_\lambda\hat\otimes{\cal H}_\lambda^*)$
in ${\cal O}(G/K\times \oline {G/K}, B(F(\lambda)))$. 
Here we used the identification $B(F(\lambda))=
F(\lambda)\otimes F(\lambda)^*$. The $U=K\times K$-spherical 
vector in $B(F(\lambda))$ is the trace in this case. 
The function $\|D_\lambda\|$ is given by 
$$\|D_\lambda\|\: G/K\to \R^+, \ \ gK\mapsto \tr( 
J_\lambda(g,0)^{-1}(J_\lambda(g,0)^{-1})^*)\ .$$

\par If $G$ is a locally compact group, then we write $\mu_G$ for 
a left invariant Haar measure on $G$. 
 
\Lemma 2.3. Let $(\pi_\lambda, {\cal H}_\lambda)$
be a unitary highest weight representation of $G$. Then 
$D_\lambda$ is square integrable if and only if 
$$\la 2\lambda+\rho,\beta\ra<0$$
for $\beta$ the highest root in $\Delta^+$. 

\Proof. That $D_\lambda$ is square integrable means
that 
$$\int_{G/Z} \|J_\lambda(g,0)^{-1}(J_\lambda(g,0)^{-1})^*)\|^2
\ d\mu_{G/Z}(gZ)<\infty.\leqno(2.4)$$
Assume for a moment that $\lambda\in i\z(\k)^*$, i.e., 
$F(\lambda)\simeq \C$. Then $J_\lambda(g,0)^2=J_{2\lambda}(g,0)$
and (2.4) becomes
$$\int_{G/Z} \|J_{2\lambda}(g^{-1},0)J_{2\lambda}(g^{-1},0)^*\|
\ d\mu_{G/Z}(gZ)<\infty.\leqno(2.5)$$
{}From the equivalence  (2)$\iff$(4) in Lemma 2.2 we obtain that (2.5) is 
equivalent to 
$$\la 2\lambda+\rho,\beta\ra<0$$
for $\beta$ the highest root in $\Delta^+$. This proves the 
lemma for the case of $\lambda\in i\z(\k)^*$. The general case 
can be easily reduced to the special case discussed above (see (2.6) below). 
\qed

\Lemma 2.4. If the highest weight $\lambda$ of 
$(\pi_\lambda, {\cal H}_\lambda)$ satisfies the 
condition of $\la 2\lambda+\rho,\beta\ra<0$, 
then $\lambda$ is regular. 

\Proof. This follows from the classification of 
unitary highest weight modules (cf.\ [EHW83]). In the scalar 
case this is proved in Lemma A.2 (iv) in Appendix A.\qed 

Assume now that $\la 2\lambda+\rho,\beta\ra<0$ holds. In view 
of Lemma 2.3  and Lemma 2.4, this implies in particular
that $D_\lambda$ is square integrable. In view 
of Lemma 2.1 we have thus proved:

\Theorem 2.5. Let $(\pi_\lambda,{\cal H}_\lambda)$ be 
a unitary highest weight representation with 
highest weight $\lambda$ satisfying the condition 
$\la 2\lambda+\rho, \beta\ra<0$ with $\beta$ the highest 
root. Then we have an onto $G$-equivariant partial isometry 

$$U_\lambda\: {\cal H}_\lambda\hat\otimes{\cal H}_\lambda^* \to 
L^2(G/K)$$
which is an isomorphism when $\lambda\in i\z(\k)^*$.
In particular, for $\lambda\in i\z(\k)^*$ 
the branching of the irreducible unitary 
representation $\pi_\lambda\otimes\pi_\lambda^*$ of 
$G\times G$ to the diagonal subgroup $G$ is completely 
continuous and multiplicity free.\qed

\subheadline{The square integrability of $\|D_\lambda\|$: general case}

In this subsection  we will prove 
a theorem that will give us the exact range of parameters for 
the square integrability of $D_\lambda$. 
For that  recall some results from [\'O{\O}88]
applied to the symmetric Lie algebra $(\s, \theta\tau)$, 
which is also compactly causal. The proof of the statements
uses simple $\s\frak{u}(1,1)$-reduction (cf.\  Appendix B). Choose the
Cartan subspace
$\t$ such that $\e\:=\t\cap \q_*$ is maximal abelian in $\q_*$.
Write $\Delta_\e=\Delta(\s_\C, \e_\C)$ for the restricted 
root system with respect to $\e_\C$. Note that $\Delta_\e$
is always of type $C_n$ or $BC_n$. Define a positive system 
$\Delta_\e^+$ by $\Delta_\e^+\:=\Delta^+\res_\e\bs\{0\}$. Further 
set $\Delta_{\e,n}^+\:=\Delta_n^+\res_\e$.
Note that $\gamma\res_{\e}\not= 0$ for
all $\gamma \in \Delta_n$. 
\par  Let $\gamma_1,\ldots ,\gamma_r$ be a  maximal set of long strongly
orthogonal roots in $\Delta_{\e,n}^+$. Then we can choose
$E^{\pm j}\in \s^{\pm \gamma_j}_{\C}$ such that with
$Y^j=-i(E^{j}-E^{-j})$ the space
$\a=\bigoplus_{j=1}^r\R Y^j$ is maximal abelian in $\p$. Let
$H^j\in i\e$ be such that $\gamma_i(H^j)=2\delta_{ij}$. Then
$$\kappa (\exp(\sum_{j=1}^rt_jY^j))=
\exp\left(-\frac{1}{2}\sum_{j=1}^r\log (\cosh (2t_j))
H^j)\right)\, .$$
Hence
$$\|D_\lambda(\exp(\sum_{j=1}^rt_jY^j))\|=
\exp \left(\frac{1}{2}\la \lambda
,\sum_{j=1}^r\log(\cosh(2t_j))H^j\ra\right)
\, .$$
Thus with $a=\exp(\sum_{j=1}^rt_jY^j)$ and $k_1, k_2\in K$:
$$\|D_\lambda(k_1ak_2)\|^2=\prod_{j=1}^r\cosh(2t_j)^{\la\lambda,H^j\ra}\, .\leqno(2.6)
$$
The root vectors $E^j$ determine a Cayley transform $c$ such that
$c(\sum_{j=1}^r\R H^j)=\a$. Recall that  $\Sigma=\Sigma (\g,\a)$  is the
set
of restricted roots of $\a$ in $\g$. Choose $\Sigma^+$ such
that $\Sigma^+\subseteq \Delta^+\circ c^{-1}$. Let
$\rho_G$ be the half sum of positive roots counted with multiplicities
and let $\rho_G^c=\rho_G\circ c^{-1}\in i\e^*$.

\Theorem 2.6. We have $D_\lambda \in L^2(G/K)$ if and only if
$$\la \lambda +\rho_G^c,\beta \ra <0$$
for $\beta$ the highest root in $\Delta^+$.

\Proof. First note that $\la \lambda +\rho_G^c,\beta\ra <0$ is
equivalent to
$$\la \lambda +\rho_G^c,\alpha\ra <0 \quad (\forall \alpha\in
\Delta_n^+)\ .$$
Further the fact that $\lambda, \rho_G^c\in i\e^*$ implies
that this is equivalent to
$$\la \lambda +\rho_G^c,\alpha\ra <0
\quad (\forall \alpha\in \Delta_{\e,n}^+)\ .$$
\par Let
$$\Delta_G(k_1ak_2):=\prod_{\alpha\in\Sigma^+}
\sinh(\alpha (\log a))^{m_{\alpha}}\, \quad k_1,k_2\in K,
a\in A$$
where $m_\alpha =\dim \g^\alpha$. Then for $f\in L^1(G)$:
$$\int_G f(g)\ dg=\int_K\int_{A^+}\int_Kf(k_1ak_2)\ \Delta_G(a)
\ d\mu_K(k_1)\  d\mu_A(a)\  d\mu_K( k_2)$$

Let $\varphi(t)=\frac{1}{2}(1-e^{-2t})$, $t\ge 0$. Then $\varphi$ is
increasing, $0\le \varphi (t)\le 1/2$, and $\varphi(t)=0$ if
and only if $t=0$. Furthermore $\sinh(t)=\varphi(t)e^{t}$.
Define
$$\Phi(a):=\prod_{\alpha\in \Sigma^+}\varphi(\alpha (\log a))^{m_\alpha}\, a\in A^+\, .$$
Then there exists a positive constant $C $ such that
$$\Phi(a)a^{2\rho_G}=\Delta (a)\le  Ca^{2\rho_G}$$
for all $a\in A^+$. Notice that by (2.6) there are constants $C_1,C_2>0$ such that
$$C_1e^{\lambda (\sum_{j=1}^rt_jH^j)}\le D_{\lambda}(k_1\exp(\sum_{j=1}^r t_jY^j)k_2)
\le C_2 e^{\lambda (\sum_{j=1}^rt_jH^j)}\, .$$
Hence for $a=\exp(\sum_{j=1}^rt_jY^j)\in A^+$
$$ C_1^2 \Phi(a) \prod_{j=1}^r\exp(2t_j\la \lambda +\rho_G^c,H^j\ra)
\leq \|D_\lambda(k_1ak_2)\|^2\Delta(a)\leq C C_2^2 \prod_{j=1}^r\exp(2t_j\la\lambda +
\rho_G^c,H^j\ra)\,.$$
The claim now follows as $t_j>0$ for $a\in A^+$.
\qed

{}From Theorem 2.6 and Lemma 2.1 we now obtain the following result:

\Theorem 2.7. Let $(\pi_\lambda,{\cal H}_\lambda)$ be 
a $G$-spherical unitary highest weight representation of $S$. 
Suppose the highest weight $\lambda$ is regular and  satisfies the
condition 

$$\la \lambda+\rho_G^c, \beta\ra<0$$
for $\beta$ the highest root in $\Delta^+$. 
Then we have an onto $G$-equivariant partial isometry

$$U_\lambda\: {\cal H}_\lambda \to L^2(G/K)$$
which is an isomorphism when $\lambda\in i(\z(\uu)\cap\q_*)^*$.
In particular, for $\lambda\in i(\z(\uu)\cap\q_*)^*$ 
the representation $\pi_\lambda\res_G$ of 
$G$ is multiplicity free with completely continuous
spectrum. \qed

\Remark 2.8. (a) The condition in Theorem 2.7 is the same as 
the one in Lemma 2.3 in the group case. The condition on the 
regularity of $\lambda$ in Theorem 2.7 is probably superfluous, 
but we did not check it. 

\par\nin (b) The results in [GrKo01] for 
$\s=\sL(2,\R)\oplus\sL(2,\R)$ show that the condition 
$\la \lambda+\rho_G^c,\beta\ra<0$
on $\lambda$ in Theorem 2.7
is also necessary for continuous branching. It would be interesting 
to know if this is generally true.\qed

\Example 2.9. Let us consider $S=\tilde \Sl(2,\R)$ with 
$G=\SO(1,1)$. In this case $\rho_G^c=0$ and this means that 
all unitary highest weight representations $\pi_\lambda$ 
have continuous branching for $\pi_\lambda\res_G$. 
\par On the other hand consider the group case 
$S=\tilde \Sl (2,\R)\times\tilde \Sl (2,\R)$. Here the condition 
$\la \lambda+\rho_G^c,\beta\ra<0$ gives a real restriction. 
Further by the results in [GrK01] we do not have continuous
branching for $\pi_\lambda\otimes\pi_\lambda^*$ to $G$ on 
the full half line. \qed

\subheadline{Applications to the classical Hardy space}

Write $\partial_s{\cal D}$ for the Shilov boundary of ${\cal D}\simeq S/U$
and recall  that $\partial_s{\cal D}$ is the  $U$-orbit 
through a certain  $E\in\partial_s{\cal D}$.  Note that ${\cal D}$ is 
a circular domain. In particular we have 
$r\partial_s{\cal D}\subeq {\cal D}$ for all $0\leq r<1$. 

\par The {\it Hardy space parameter} $\lambda_h\in i
(\z(\uu)\cap \q_*)^*$
is defined by 
$$\lambda_h=-\rho_n\:=-{1\over 2}\sum_{\alpha\in \Delta_n^+}\alpha\ .$$
The corresponding Hilbert space ${\cal H}_{\lambda_h}$ is the 
{\it classical Hardy space} on ${\cal D}$: 
$$ {\cal H}^2({\cal D})\:=\{ f\in {\cal O}({\cal D})\: 
\sup_{0\leq r<1} \int_U |f(ru(E))|^2\ d\mu_U(u) <\infty\}\ .\leqno(2.7)$$
Here we identified the Shilov boundary of $\partial_s{\cal D}$
with $U(E)$ as explained before. 
We write $\pi_h\:=\pi_{\lambda_h}$ for the unitary 
irreducible representation of $S$ on ${\cal H}^2({\cal D})$. 

\Lemma 2.10. Assume that $(\s,\g,\h)$ is a causally symmetric 
triple. Then $\lambda_h$ is regular
and satisfies the condition
$\la \lambda_h+\rho_G^c, \beta\ra <0$.

\Proof. See Proposition A.3 in the appendix.
\qed

{}From this  we obtain the following important corollary to 
Theorem 2.8: 

\Corollary 2.11. Let $(\s,\g,\h)$ be a causally symmetric triple
and ${\cal H}^2({\cal D})$ be the classical
Hardy space on ${\cal D}\simeq S/U$. Then the inverse of the 
Segal-Bargmann transform 
$$U\: (\pi_h\res_G, {\cal H}^2({\cal D})) \to (L,  L^2(G/K))$$
is a $G$-equivariant unitary isomorphism. 
In particular, the branching $\pi_h\res_G$ 
is completely continuous and multiplicity free.\qed

\sectionheadline{3. Analytical and geometrical constructions
on $\Xi$}

{}From now on we will assume that $(\s,\g,\h)$ is a causally
symmetric triple. Further we will require $S_\C$ to be simply connected
and henceforth $S$,$U$,$G$,$K$ and $H$
will denote the analytic subgroups of $S_\C$ with Lie algebras 
$\s$, $\uu$, $\g$, $\k$ and $\h$.

\par We have seen in Section 1 that $\partial_s\Xi\simeq G/H$ is open and 
dense in $\partial_s {\cal D}$ (cf. \ Theorem 1.14). Now on $\partial_s\Xi$ 
we have a natural $G$-invariant measure $\mu_{G/H}$ while on $\partial_s{\cal D}$
we have a natural $U$-invariant measure $\mu_{\partial_s{\cal D}}$. These two 
measures are related through a density function, i.e. $d\mu_{G/H}(z)=
{1\over |\psi (z)|^2}d\mu_{\partial_s{\cal D}}(z)$. 
In this section we will show that one can choose the measurable function $\psi$ on $\partial_s{\cal D}$
in such a way that it admits a continuous extension to ${\cal D}\amalg  \partial_s{\cal D}$ which in addition 
is holomorphic on ${\cal D}$. The function $\psi$ will satisy 
a cocycle property  which allows us to 
identify the $L^2$-spaces $L^2(G/H)$ and $L^2( \partial_s{\cal D}, \mu_{\partial_s{\cal D}})$
in a natural $G$-equivariant way (cf.\ Lemma 3.11 below). Moreover the fact that 
$\psi$ has a holomorphic extension to $\Xi\simeq {\cal D}$ will be used for
comparing the classical Hardy 
space from (2.7) with our Hardy space on $\Xi$ in Section 5.

\subheadline{The function $\psi$}

In this subsection we give the construction of the function $\psi$ 
mentioned above. For that we first have to collect some facts on certain  
finite dimensional $K_\C$-spherical representations of $U_\C$.

\par Denote by $(\pi_m,V_m)$ the irreducible representation
of $U_\C$ with lowest weight $-m\rho_n$. If such a representation
exists it is $K_\C$-spherical. Note that $V_m\simeq \C$ is one 
dimensional since $\rho_n\in i(\z(\uu)\cap\q_*)^*$. 
Also notice that the representation
$\det \Ad|_{\p^-}$ of $U$  has lowest weight
$-2\rho_n$. Let
$$\zeta =\frac{1}{2}(\gamma_1+\ldots +\gamma_r)\in i(\z(\uu)\cap\q_*)^*$$
where $\gamma_1,\ldots ,\gamma_r$ is a maximal
set of long strongly orthogonal roots in $\Delta_n^+$. Let
$$d=\dim \s_\C^{\frac{1}{2}(\gamma_i+\gamma_j)}
=\dim \s_\C^{\frac{1}{2}(\gamma_i-\gamma_j)}\, $$
(cf.\ Table II in Appendix A). Then
$$\rho_n=(1+\frac{d(r-1)}{2})\zeta \leqno(3.1)$$
(cf. Appendix A).

\Proposition 3.1. Let the notation be as above. Let 
$m\in\N$ be minimal such that $(\pi_m, V_m)$ exists. Then
the following holds:
\item{(i)} If $(\g,\h)$ is of Cayley type and
$\g \not= \so(2,2k+1),\sp(2n,\R)$ then
$m=1$.
\item{(ii)}  If $(\g,\h)$ is of Cayley type and
$\g = \so(2,2k+1),\sp(2n,\R)$ then
$m=2$.
\item{(iii)} If $(\s,\g)=(\su(2n,2n),\sp(n,n))$ then
$m=1$.
\item{(iv)}  If $(\s,\g)=(\sp(2n,\R),\sp(n,\C))$ then
$m=2$.
\item{(v)} If $(\s,\g)=(\so(2,n),\so(1,n))$ for $n>2$ then
$m=1$ if $n$ is even and $m=2$ for $n$ odd.

\Proof. Using Table II for $d$ in Appendix A, this follows from
(3.1) by straightforward computation. \qed

\Lemma 3.2. The set $G_\C U_\C P^-$ is open and dense in
$S_\C$. Furthermore $S\subseteq G_\C U_\C P^-$.

\Proof. This is Theorem 2.4 in [\'O{\O}88].
\qed

Assume now that $(\pi_m,V_m)$ exists. Then $V_m=\C v_m$ and
$v_m$ is both a lowest weight vector and a $K_\C$-spherical
vector. Normalize $v_m$ such that  $(v_m,v_m)=1$.
Let
$${\cal U}:=\{z\in \p^+\:\exp(z)\in G_\C U_\C P^-\}\,
.$$
Since $G_\C U_\C P^-$ is open and dense in $S_\C$, it follows
from Lemma 3.2 that ${\cal U}$ is an
open dense and locally $G_\C$-invariant subset of $\p^+$.
Moreover  ${\cal D}\subeq {\cal U}$.

For $s\in {\cal U}$ define
$K_\C u_G(s)\in K_\C \backslash U_\C$ by
$$ s\in G_\C u_G(s)P^-\, .\leqno(3.2)$$
Define a function $\psi_m: {\cal U}\to \C$ by
$$\psi_m(z):=(\pi_m(u_G(\exp(z))v_m,v_m)\, .\leqno(3.3)$$

\Lemma 3.3. The function $\psi_4$ extends
to an holomorphic function on $\p^+$.

\Proof. This follows from  [B\'O01], Proposition 4.1 and the
table on p. 294, but we will give a short proof here. First
notice that there exists a finite
dimensional representation $(\sigma,W )$ of $S_\C$ with
lowest weight $-2\rho_n$. As $-2\rho_n\in i(\z (\frak{u} )\cap \q_*)^*$
it follows by Helgason's Theorem that the irreducible representation
with lowest weight $-4\rho_n$ is $G_\C$-spherical. Denote this
representation by $(\sigma_4,W_4)$. Let $v\in W_4$ be a lowest weight vector
and $u\in W_4$ be a $G_\C$-invariant vector. We can normalize $v$ and $u$ such
that $(u,v)=1$. It then follows that
$$(\forall s\in S_\C)\qquad (\sigma_4(s)v,u)=\left(\sigma_4\left(u_G\left(s
\right)\right)v,u\right)\, .$$
Now note that the representation of $U_\C$ generated by
$\sigma_4(U_\C)v$ is in fact irreducible and hence
equivalent to $\pi_4$. Thus
$$(\forall z\in {\cal U})\qquad (\sigma_4\left(u_G(\exp(z))\right)v,u)=\psi_4(z)$$
and the claim follows as $z\mapsto (\sigma_4(\exp z)v,u)$ is holomorphic on
$\p^+$.\qed

\Proposition 3.4. Let $\psi_m$ be as above. Then the following
assertions hold:
\item{(i)} Let $g\in G_\C$ and $z\in {\cal U}$ such that $gz\in {\cal
U}$.
Then
$$\psi_m(g z)=J_{m\rho_n}(g,z)\psi_m(z)\, .$$
\item{(ii)} $\psi_m$ is holomorphic and has no zeros on ${\cal D}$.

\Proof. (i) We have $g\exp (z)=\exp (g z)J(g,z)p$ with
$p\in P^-$. Hence
$$\exp (g z)= g\exp (z)J(g,z)^{-1}p^\prime$$
for some $p^{\prime} \in P^-$. It follows that
$$\psi_m(gz)=(\pi_m(\exp (z))\pi_m(J(g,z)^{-1})v_m,v
_m))
= J_{m\rho_n}(g,z)\psi_m(z)$$
because $v_m$ is a weight vector with weight $-m\rho_n$.
\par\nin (ii) The holomorphicity of $\psi_m$ is clear by 
construction. To see that $\psi_m(z)\neq 0$ for $z\in {\cal D}$ write 
$z=g(0)$ for some $g\in G\exp(i\Omega)$. Since ${\cal D}\subeq 
{\cal U}$, we obtain from (i) that 

$$\psi_m(z)=\psi_m(g(0))=J_{m\rho_n}(g,0)\psi_m(0)\ .$$ 
Now $\psi_m(0)=(v_m, v_m)=1\neq 0$ and $J_{m\rho_n}(g,0)\neq 0$
by construction. \qed

{}For the rest of this subsection we consider $\psi_m$ as a function on ${\cal D}$ only.
Note that $\psi_2$ always exists by Proposition 3.1.
Since ${\cal D}$ is simply connected and $\psi_m$ is zero-free
it follows from Proposition 3.4(ii) that we can define
a holomorphic square root $\psi(z)\:=\sqrt{\psi_2(z)}$ which
becomes unique under the requirement $\psi(0)=1$.
Note that $\psi=\psi_1$ in case $(\pi_1, V_1)$ exists.

\par Write $\tilde S$ for the universal  covering of $S$ and
$G_1$ for the analytic subgroup of $\tilde S$ with Lie algebra $\g$.
Note that $J_{\rho_n}(s,z)$ exists for all $s\in \tilde S$
and $z\in {\cal D}$ (cf.\ Section 2). From Proposition 3.4(i)
we thus obtain that $\psi(gz)=J_{\rho_n}(g,z)\psi(z)$ for all
$g\in G_1$ and $z\in {\cal D}$. In particular,
$J_{\rho_n}(g,z)={\psi(z)\over \psi(gz)}$ and so
$J_{\rho_n}(g,z)$, initially defined only on $G_1\times {\cal D}$,
factors to a function on $G\times {\cal D}$, which
we denote by $J_h^{-1}$.
Summarizing our discussions we have proved:

\Proposition 3.5. There exist an unique holomorphic function
$\psi\:{\cal D}\to \C^*$ and a unique  analytic function
$J_h\: G\times {\cal D}\to \C^*$ with the following properties:
\item{(i)} $\psi(gz)=J_h^{-1}(g,z) \psi(z)$ for all $g\in G$,
$z\in {\cal D}$.
\item{(ii)} $\psi^2=\psi_2$ and $\psi(0)=1$.
\item{(iii)} $J_h$ satisfies the cocycle property $J_h(g_1g_2,z)=
J_h(g_1, g_2z)J_h(g_2,z)$ for $g_1,g_2\in G$, $z\in {\cal D}$ as well
as $J_h^{-2}=J_{2\rho_n}$ and $J_h(\1,z)=1$ for all
$z\in {\cal D}$. \qed

\Example 3.6. Let $S=G\times G$, and $G$ the diagonal
group. Write $G\subset P_G^+K_\C P_G^-$ for the
triangular decomposition for $G_\C$ and let $k_1(g)\in K_\C$ be
the corresponding Harish-Chandra projection.
If ${\cal D}_G\subeq \p_G^+$ is the bounded realization of
$G/K$, then ${\cal D}={\cal D}_G\times {\cal D}_G^{\rm opp}$
where ${\cal D}_G^{\rm opp}$ denotes ${\cal D}_G$ equipped
with the opposite complex structure. Note that $P^+=P_G^+\times P_G^-$
and $P^-=P_G^-\times P_G^+$. In particular we are realizing
$S/U$ inside $(G_\C /K_\C P^-) \times (G_\C /K_\C P^+)$. Recall
that the conjugation $g\mapsto \oline g$ with respect to $G$
satisfies $\oline{P^+}=P^-$. 
\par Let $(s,t)\in G\times G$. Then
$$(s,t)=(t,t)(t^{-1}s,\1)=(t,t)(p_1^+k_1(t^{-1}s)p_1^-,\1)
=(tp_1^+,tp_1^+)(k_1(t^{-1}s)p_1^-,(p_1^+)^{-1})
\, .$$
Hence 
$$ G_\C u_G(s,t)=G_\C (k_1(t^{-1}s),\1)\ .$$ 
As $U_\C = K_\C\times K_\C$ it follows that any  $K_\C$-spherical
representation $\pi_U$ of $U_\C$ is of the form $\pi_U =\pi\otimes
\pi^*$
where $\pi$ is an irreducible representation of $K_\C$. Hence
$$\psi(z,w)=\pi_1(k_1(\exp(-\bar{w})\exp (z)))=K_{-\rho_n^G }(z,w)^{-1}\qquad 
(z,w\in {\cal D}_G)\leqno (3.4)$$
where $\rho_n^G$ is the $\rho_n$ for the Hermitian group
$G$. Notice that $K_{-\rho_n^G }$ is nothing else but the 
reproducing kernel of the classical Hardy space on ${\cal D}_G$. \qed

\subheadline{$\SU(1,1)$-reduction}

In this subsection we will discuss the
case $S=\SU(1,1)$ even if it is not in
our list of causally symmetric triples. The reason
is, that many calculations can be reduced to this
situation using  $\tau$-equivariant
embedding of $\SU(1,1)$ into $S$. Let
$$S =\SU(1,1)=\{\pmatrix{\alpha & \beta\cr
\bar{\beta} & \bar{\alpha}\cr}\: \alpha, \beta \in \C,\,
|\alpha|^2-|\beta|^2=1\}$$
Define $\tau :S\to S$ by conjugation with  $\pmatrix{0 & 1\cr 1 & 0\cr}$.
On $S_\C = SL(2,\C)$ the involution $\tau$ is given by
$$\tau\pmatrix{a & b\cr c & d\cr}=
\pmatrix{ d & c\cr b & a\cr}$$gko6.
In particular
$$G=\{\pmatrix{ \cosh(t) & \sinh(t)\cr
\sinh(t)& \cosh(t)\cr}\: t\in \R\}\, .$$
In this situation we have
$$Y^0=\frac{1}{2}\pmatrix{0 & 1\cr 1 & 0\cr}, \quad 
X^0=-i\frac{1}{2}\pmatrix{1 & 0\cr 0& -1\cr}\quad \hbox{and}\quad 
Z^0= \frac{1}{2}\pmatrix{0 & -i\cr i & 0\cr}. $$

\nin For the general situation of a causally symmetric triple recall 
the $\su(1,1)$-triple $\{X^j, Y^j, Z^j\}$ from Appendix B. 
{}From the relations (B.1) and (B.2) in Appendix B
it now follows:

\Lemma 3.7. Let the notation be as above. Then
the map $X^j\mapsto \pmatrix{1 & 0\cr 0 & -1\cr}$,
$Y^j\mapsto \pmatrix{0 & 1\cr 1 & 0\cr}$,
and $Z^j\mapsto \pmatrix{0 & -i\cr i & 0\cr}$
defines a Lie algebra homomorphism into
$\s\frak{u} (1,1)$ intertwining the involution $\tau$ on
$\s$ and the above involution on $\s\frak{u} (1,1)$. Furthermore
it also intertwines the Cartan involution $\theta$ on
$\s$ and the Cartan involution $X\mapsto -X^*$
on $\su (1,1)$.
\qed

On the group level we have:
$$U_\C=\{\pmatrix{\gamma & 0\cr 0 & \gamma^{-1}\cr}\: \gamma\in \C^*\}\ ,$$
$$ P^+=\{\pmatrix{1 & z\cr 0 & 1\cr}\: z\in \C\}\ ,$$
and
$$P^-=\{\pmatrix{1 & 0\cr w & 1 \cr}\: w\in \C\}\, .$$
Furthermore
$$\pmatrix{1 & z\cr 0 & 1\cr}\pmatrix{\gamma & 0 \cr 0 & \gamma^{-1}\cr}
\pmatrix{1 & 0\cr w & 1\cr}=
\pmatrix{\gamma +zw\gamma^{-1} & z\gamma^{-1}\cr
\gamma^{-1}w & \gamma^{-1}\cr}\, .$$
Consequently, if $d\not= 0$ we have
$$\pmatrix{a & b\cr c & d\cr}
=\pmatrix{1 & b/d\cr 0 & 1\cr}\pmatrix{\frac{1}{d} &  0\cr
0 & d\cr}\pmatrix{1 & 0\cr c/d &1 \cr}\, .$$
For the $G_\C U_\C P^-$-decomposition we notice that
$$\pmatrix{a & b\cr b & a\cr}\pmatrix{\gamma & 0 \cr 0 & \gamma^{-1}\cr}
\pmatrix{1 & 0\cr w & 1\cr}=
\pmatrix{a\gamma +bw\gamma^{-1} & b\gamma^{-1}\cr
b\gamma+  a\gamma^{-1}w & a\gamma^{-1}\cr}\, .$$
Hence, using  that $a^2-b^2=1$, we obtain that 
$$\pmatrix{a & b\cr c & d\cr}=
\pmatrix{ \frac{d}{\sqrt{d^2-b^2}}& \frac{b}{\sqrt{d^2-b^2}}
\cr \frac{b}{\sqrt{d^2-b^2}} & \frac{d}{\sqrt{d^2-b^2}} \cr}
\pmatrix{\frac{1}{\sqrt{d^2-b^2}} & 0
\cr 0 & \sqrt{d^2-b^2}\cr}
\pmatrix{1 & 0\cr \frac{cd-ab}{d^2-b^2} & 1\cr}\, .$$
Thus the open sets $G_\C K_\C P^-$ and $P^+K_\C P^-$ in $S_\C$ 
are given by:
$$G_\C U_\C P^- =\{\pmatrix{a & b\cr c & d\cr}\: b^2-d^2\not= 0\} $$
$$P^+U_\C P^- = \{\pmatrix{a & b\cr c & d\cr}\mid d\not= 0\}$$
Furthermore we notice that the $P^+, U_\C$, and
$P^-$ components,
whenever defined, are
given by
$$s=\pmatrix{a & b\cr c & d\cr}\mapsto u_G(s)=
\pmatrix{\frac{1}{\sqrt{d^2-b^2}}& 0 \cr 0&\sqrt{d^2-b^2}\cr} $$
and
$$s=\pmatrix{a & b\cr c & d\cr}\mapsto\pmatrix{1 & b/d\cr 0 & 1\cr}
\pmatrix{d^{-1}& 0 \cr 0& d\cr}
\pmatrix{1 & 0\cr c/d& 1\cr}\, .$$
Notice that a double covering is needed in general for $u_G$, but
as an element of the coset $K_\C \bs U_\C$, with $K_\C =\{\pm 1\}$
it is well defined. In the following we will identify 
$P^+$ with $\C$ by $z\mapsto \pmatrix{1 & z\cr 0 &1\cr}$ and
$\C^*$ with $U_\C$ by
$\gamma \mapsto \pmatrix{\gamma & 0 \cr 0 & \gamma^{-1}\cr}$. 
In particular we get the following lemma:

\Lemma 3.8. Identify $\p^+$ and $\p^-$ with $\C$, and
$U_\C$ with $\C^*$ in the way explained above. Let $z,w\in \C$ and
$g=\pmatrix{a & b\cr c & d\cr}$. Then the following assertions hold:
\item{(i)} If $z\bar{w}\not=1$ then ${\cal K}(z,w)\in U_\C$ is
defined and given by
$${\cal K}(z, w)=\frac{1}{1-z\bar{w}}\, .$$
\item{(ii)} If $az+d\not= 0$ then $J(g,z)$ is defined and
$$J(g,z)=cz+d\ .$$
\item{(iii)} Assume that $|z|<1$. Then
$u_G(z)$ is defined and
$$u_G(z)=\frac{1}{\sqrt{1-z^2}}\, .$$
\qed

\subheadline{ Description of $\partial_s\Xi$ in $\partial_s{\cal D}$}

In this subsection we will use the $\su(1,1)$-reduction 
to identify $\partial_s\Xi$ in $\partial_s{\cal D}$ as 
the non-vanishing locus of the function $\psi$.  Recall
that $\psi_4$ extends to a holomorphic function on
$\p^+$.

\Lemma 3.9. We have:
$$\partial_s\Xi=\{z\in \partial_s{\cal D}\:\psi_4(z)\not= 0\}\ . $$

\Proof. ``$\subeq$'': Recall the element $z_1\in \partial_s\Xi$
with isotropy subgroup $H$, i.e., $\partial_s\Xi=G(z_1)\simeq 
G/H$. We first show that $\psi_4(z_1)\neq 0$
by using $\SU(1,1)$-reduction. For $S=\SU(1,1)$ we have 
${\cal D}=\{z\in\C\: |z|<1\}$ and $z_1=i$. In particular 
it follows from Lemma 3.8(iii) that $\psi_4(z_1)\neq 0$.
Now $\psi_4(z_1)\neq 0$ in 
the general case follows from simple $\SU(1,1)$-reduction (cf.\ Lemma 3.7 
and our structural results in Appendix B).  
{}From the covariance property of $\psi_4$ (cf.\ Proposition 3.5)
we hence get that $\psi_4(gz_1)\neq 0$ for all $g\in G$. In particular 
$\psi_4$ has no zeros on $G(z_1)=\partial_s\Xi$.

\par\nin  ``$\supeq$'': 
By Theorem 1.14 we know that $\partial_s\Xi$ is open and dense 
in $\partial_s{\cal D}$. 
Let  $Y^j$ and $\a = \sum_{j=1}^r\R Y^j$ be as in Appendix B.
Let $A=\exp (\a )$. Then $G=KAH$.
Let $g_n$ be a sequence in $G$ such
that $g_n(z_1)\to \partial_s{\cal D}\bs \partial_s\Xi$.
Write $g_n=k_{n}a_nh_{n}$ with
$k_{n}\in K,\, h_n\in H$ and $a_n=\exp (\sum_{j=1}^r t_{r,n}Y^j)\in A$.
Note that $g_n(z_1)\to \partial_s{\cal D}\bs \partial_s\Xi$
means precisely that $a_n\to \infty$. 
\par We have $g_n(z_1)=k_na_n(z_1)$. We also have that
$|J_{4\rho_n}(k,z)|=1$ for all $k\in K$.
Hence it  follows from Proposition
3.5 that
$$|\psi_4(k_na_n(z_1))|=|J_{4\rho_n}(k_n,a_n(z_1))|\cdot
|\psi_4(a_n(z_1))|
=|\psi_4(a_n(z_1))|= c|J_{4\rho_n}(a_n,z_1)|$$
with $c\:=|\psi_4(z_1)|>0$.
By  $\SU(1,1)$-reduction (cf.\ Lemma 3.7 and Lemma 3.8(ii)) we get
$$J_{4\rho_n}(a_n,z_1)=\prod_{j=1}^r
(i\sinh(t_{n,j})+\cosh(t_{n,j}))^{-4\rho_n(H^j)}\, .$$
As $a_n\to\infty$,  there
exists an index $j$ and a subsequence
$s_n=t_{n,j}$ such that $|s_n|\to \infty$. Therefore
$$\prod_{j=1}^r\left(i\sinh(t_{n,j})+\cosh(t_{n
,j})\right)^{-4\rho_n(H^j)}\to
0$$
because $\rho_n(H^j)>0$. This concludes the proof of the Lemma. \qed

\Theorem 3.10. The function $\psi$ extends to a continuous
function on ${\cal D}\cup \partial_s{\cal D}$
such that
$$\partial_s\Xi =\{z\in\partial_s{\cal D}\: \psi(z)\not=0\}\, .$$

\Proof. Notice that the set ${\cal D}\cup \partial_s\Xi$ is simply connected. Thus Lemma 3.9 
implies that there exists a unique continuous function $\psi: {\cal D}\cup\partial_s\Xi\to \C^*$ such
that $\psi(0)=1$ and $\psi(z)^4=\psi_4(z)$. In particular this function
agrees with our old definition of $\psi$ on $\cal D$. For $z\in \partial_s{\cal D}
\setminus \partial_s\Xi$ define
$\psi(z)=0$. It follows from Lemma 3.9 that $\psi$ is continuous.\qed

\subheadline{$L^2$-isomorphism on the Shilov boundaries}

Let $\mu_{\partial_s{\cal D}}$
be the unique (up to constant) $U$-invariant
measure on $\partial_s{\cal D}$. If $\lambda_h=-\rho_n$ is
analytically integral, then 
$S$ acts unitarily on $L^2(\partial_s{\cal D},\mu_{\partial_s{\cal
D}})$
by
$$(\pi_h(s)f)(z)=J_h (s^{-1},z)^{-1}f(s^{-1}z)\qquad (s\in S, 
\ z\in {\cal D},\  f\in {\cal H}^2({\cal D})) \leqno(3.5)$$
with $J_h =J_{-\rho_n}$. 
In the case where $m=2$ (see Proposition 3.1) we need
to go to a double covering of $S$. But notice, that $\pi_h(s)$
is always defined for $s\in G$ according to Proposition 3.5. 
One of the consequences of this observation is the following
identification of $L^2$-spaces:

\Lemma 3.11. The following assertions hold:
\item{(i)} The $G$-invariant measure on $G/H\simeq \partial_s\Xi$ as a
subset of $\partial_s{\cal D}$ is given by
$$f\mapsto \int_{\partial_s{\cal D}}f(z)|\psi(z)|^{-2}\, 
d\mu_{\partial_s{\cal D}} (z)$$
where $d\mu_{\partial_s{\cal D}}$ is the unique (up to constant)
$U$-invariant
measure on $\partial_s{\cal D}$.
\item{(ii)} The mapping $f\mapsto {1\over \psi}f$ is a $G$-equivariant
isomorphism of $L^2(G/H)$ onto $L^2(\partial_s{\cal D},
\mu_{\partial_s{\cal D}})$.

\Proof. This follows from the covariance of $\psi$ in Proposition 3.5
(cf. [B\'O01], Theorem 5.1).
\qed

\sectionheadline{4. The classical Hardy space inside $L^2(G/H)_{\rm mc}$}

In this section we will show
how to realize the classical Hardy space 
${\cal H}^2({\cal D})$  (cf.\ (2.7)) in the most continuous spectrum 
of $L^2(G/H)$.  

\subheadline{The most continuous spectrum $L^2(G/H)_{\rm mc}$}

Recall the involution $\sigma$ on $\g$
associated to $\h$. Let $\g=\h\oplus\q$
be the $\sigma$-eigenspace decomposition. 
Note that we can choose $\a\subeq \q\cap\p$
since $(\g,\sigma)$ is  non-compactly causal (cf.\ [H\'O96]). 
In particular, the minimal
parabolic subgroup $P_{\rm min}=MAN$ is also a minimal 
{\it $\theta\sigma$-stable} parabolic subgroup of $G$.  Here, as usual,
$M=Z_K(\a)$. If $\m$ is the Lie algebra of $M$, then note that 
$\m\subeq \h$ since $\a\subeq \q\cap\p$. Moreover we have 
$$M=Z_H(\a)\subeq H$$
as $G_\C$ is simply connected (cf.\ [H\'O96,  Lemma 3.1.22]). 

Write $\pi_\delta=\Ind_{MAN}^G (\delta)$ for  the 
$K$-spherical unitary principal series of parameter $\delta\in i\a^*$.
Notice that $\pi_\delta \simeq \pi_{w\delta}$ for $w \in
{\cal W}$.
The {\it most 
continuous part} $L^2(G/H)_{\rm mc}$ in $L^2(G/H)$ is by definition 
the $G$-invariant subspace in $L^2(G/H)$ which corresponds in 
the Plancherel formula  to all 
principal series induced from a minimal $\theta\sigma$-stable
parabolic subgroup. Denote by ${\cal W}_0=N_{K\cap H}(\a)/ Z_{K\cap
H}(\a)$
the very little Weyl group. 
Since $M=M\cap H$ we now obtain from [BS97] or [D98]
that 

$$(L,L^2(G/H)_{\rm mc})\simeq 
\Big(\int_{i\a^*/{\cal W}}^\oplus \pi_\delta\otimes \id
\ d\mu(\delta), \int_{i\a^*/{\cal W}}^\oplus {\cal H}_\delta\otimes 
\C^{|{\cal W}/{\cal W}_0|}
\ d\mu(\delta)\Big)$$
where $\mu$ is a Borel measure that is completely continuous
with respect to the Lebesgue measure on $i\a^*$. 
In particular, up to a set of measure zero, 
the multiplicity of $\pi_\delta$ in 
$L^2(G/H)_{\rm mc}$ is $|{\cal W}/{\cal W}_0|$. 
Let $\omega\subeq i\a^*/{\cal W}$ be an open  subset. 
Then we write 

$$P_\omega\: L^2(G/H)_{\rm mc}\to \int_\omega {\cal H}_\delta\otimes 
\C^{|{\cal W}/{\cal W}_0|}
\ d\mu(\delta)$$
for the orthogonal $G$-invariant projection associated to $\omega$. 

\Definition 4.1. We call a closed $G$-invariant subspace of 
$W\subeq L^2(G/H)_{\rm mc}$ a subspace with {\it full spectrum}
if $P_\omega(W)\neq \{0\}$
for all non-empty open subets $\omega\subeq i\a^*/{\cal W}$. \qed 

\par It is interesting to compare the most-continuous 
spectrum $L^2(G/H)_{\rm mc}$ with the spectrum of  $L^2(G/K)$.
The Plancherel Theorem for $L^2(G/K)$ has the form 

$$(L,L^2(G/K))\simeq \int_{i\a^*/{\cal W}} \pi_\delta\
d\mu'(\delta)$$
with $\mu'$ a Borel measure completely continuous
with respect to $d\delta$. Hence we see that the spectrum 
of $L^2(G/K)$ is the same as the most
continuous part of the  spectrum of $L^2(G/H)$,
but the difference is that the spectrum of $L^2(G/K)$ is 
multiplicity free. 

\subheadline{Embedding of ${\cal H}^2({\cal D})$ in the most 
continuous spectrum}

Recall the $U$-invariant measure $\mu_{\partial_s{\cal D}}$
on the Shilov boundary of ${\cal D}$
and the boundary value map of the classical 
Hardy space:  

$$b\: (\pi_h, {\cal H}^2({\cal D}))\to (\pi_h, 
L^2(\partial_s{\cal D}, \mu_{\partial_s{\cal D}})), \ \ 
f\mapsto (z\mapsto \lim_{r\to 1\atop r<1} f(rz))$$
which is an $S$-equivariant isometric embedding. 
On the other hand we have the $G$-equivariant 
isomorphism from Lemma 3.11:

$$\Psi\: (\pi_h, L^2(\partial_s{\cal D}))\to (L, L^2(G/H)), 
\ \ f\mapsto \psi f\ .$$

\Theorem 4.2. The  mapping 
$$\Psi\circ b\: (\pi_h\res_G, {\cal H}^2({\cal D}))\to (L, L_{\rm
mc}^2(G/H))
$$
is a $G$-equivariant isometric embedding. Moreover 
the image $\im \Psi\circ b$ is a mulitiplicity free subspace 
of full spectrum in $L^2_{\rm mc}(G/H)$. 

\Proof. From Lemma 3.11 it is clear that $\Psi\circ b$
is a $G$-equivariant map from ${\cal H}^2({\cal D})$ into 
$L^2(G/H)$. It remains to specify the image. From 
Corollary 2.11 we obtain that $\pi_h\res_G\simeq L^2(G/K)$. Thus 
the remaining assertions of the theorem follow from 
our discussion of the most continuous spectrum in the previous
subsection. \qed

\sectionheadline{5. The Hardy space on $\Xi$}

This final section of the paper is devoted to the Hardy space ${\cal H}^2(\Xi)$ on $\Xi$. 
After a brief disgression on compression semigroups of ${\cal D}\simeq \Xi$ we 
give the definition of the Hardy space and show that it is in fact a Hilbert 
space. Also we show the existence of a boundary value map $b\: {\cal H}^2(\Xi)\to L^2(G/H)$. 
Subsequently, using the results of Section 3, we show that there is a  natural identification 
of ${\cal H}^2(\Xi)$ and ${\cal H}^2({\cal D})$. This together with Theorem 4.2 will then give us that 
$\im b \subeq L^2(G/H)_{\rm mc}$, the main result of this paper.

\subheadline{Definition of ${\cal H}^2(\Xi)$ and first properties}

In order to define the Hardy space on $\Xi$ we first have 
to recall some facts on semigroups compressing ${\cal D}$, resp. $\Xi$.

\par Recall the space $\e=\t\cap\q_*$ and associated 
restricted root system $\Delta_\e$ from Section 2. 
For $\alpha\in \Delta_\e$ we denote by $\check\alpha\in i\e$
the coroot of $\alpha$. 
Define 
$$C_{\rm min}\:=\cone(\{\check \alpha\: \alpha \in \Delta_{\e,n}^+\})\ ,$$
where $\cone(\cdot)$ refers to the convex cone generated by $(\cdot)$. 
Then 
$$W_{\rm min}\:=\Ad(G)C_{\rm min}$$
is a minimal open convex cone in $i\q_*$ and 
$$\Gamma\:=G\exp(W_{\rm min})$$
is a $G$-biinvariant subsemigroup of $S_\C$ (cf.\ [H\'O96, Ch.\ 4]). 
The closure of $\Gamma$ in $S_\C$ is given by 
$\oline \Gamma=G\exp(\oline {W_{\rm min}})$. The polar mapping 
$$G\times \oline{W_{\rm min}}\to \oline \Gamma, \ \ (g,X)\mapsto g\exp(X)$$
is a homeomorphism (cf.\ [N99, Th.\ XI.1.7]). 
Moreover, $\oline \Gamma$ is an involutive 
semigroup with involution 
$$\oline\Gamma\to\oline\Gamma,
\ \  \gamma=g\exp(X)\mapsto \gamma^*=\exp(X)g^{-1}\ .$$
Note that $\Gamma^{-1}$ compresses ${\cal D}$: 

$$(\forall \gamma\in \Gamma^{-1})\qquad \gamma(\oline {\cal D})
\subeq {\cal D}$$
(cf.\ [N99, Th.\ XII.3.3]).
{}From the realization of $\oline \Xi$ in $\oline{\cal D}$ we therefore 
obtain an action of $\Gamma^{-1}$ on $\oline\Xi$ with the 
property 

$$(\forall \gamma\in \Gamma^{-1})\qquad 
\gamma(\oline \Xi)\subeq \Xi.\leqno (5.1)$$

\Definition 5.1. {\bf (Hardy space on $\Xi$)} The Hardy 
space on $\Xi$ is defined as 
$${\cal H}^2(\Xi)=\{ f\in {\cal O}(\Xi)\: \|f\|^2\:=
\sup_{\gamma\in \Gamma}\int_{G/H} |f(\gamma^{-1} gz_1)|^2 
\ d\mu_{G/H}(gH)<\infty\}. \qeddis

{}From the definition of the Hardy space it is not clear 
yet that ${\cal H}^2(\Xi)$ is a Hilbert space. 
This will follow from the following geometric fact and 
its proof:

\Lemma 5.2. The semigroup orbit $\Gamma^{-1}(z_1)\subeq \Xi$ is 
open and measure dense. 

\Proof. Recall the elements $H^j\in i\e$ and note that 
$\sum_{j=1}^n t_j H^j\in C_{\rm min}$ whenever $t_j\geq 0$
and not all $t_j=0$. By simple $\SU(1,1)$-reduction 
we obtain for all $j$ that 

$$\exp(-t_j H^j)(z_1)=\exp(i\arctan ( e^{-t_j}) Y^j)(0).\leqno(5.2)$$
Write $D\:=\{z\in \C\: |z|<1\}$ for the unit disc
and identify $A\exp(i\Omega)(0)$ with $D^n$ in the 
obvious way through $\SU(1,1)$-reduction. 
With $D^-=D\bs \R$ equation (5.2) then gives  
$$A{\cal W}\exp(-C_{\rm min})(z_1) =(D^-)^n.\leqno(5.3)$$
In particular $A{\cal W}\exp(-C_{\rm min})(z_1)$
is  open and of full measure in $A\exp(i\Omega)(0)=D^n$. 
Sweeping out $\Xi$ by $G$-orbits through $\exp(i\Omega)(0)$
proves the lemma.\qed 

\par If $M$ is a complex manifold,  then we write 
${\cal O}(M)$ for the Fr\'echet space of holomorphic functions
on $M$.

\Corollary 5.3. The Hardy space ${\cal H}^2(\Xi)$ is
a Hilbert space of holomorphic functions. Moreover,
the inclusion mapping ${\cal H}^2(\Xi)\into {\cal O}(\Xi)$
is continuous. 

\Proof. If we know that ${\cal H}^2(\Xi)$ is complete, 
then the standard semigroup techniques imply that 
${\cal H}^2(\Xi)$ is a Hilbert space (cf.\ [H\'O\O91] or [N99, Ch.\ XIV]). 
The completeness and the continuity of the embedding 
${\cal H}^2(\Xi)\into {\cal O}(\Xi)$  is easily obtained 
from (5.3) and Cauchy's Theorem (use $\Gamma=G\exp(C_{\rm min}) G$ and choose 
local coordinates -- recall the definition of the Hardy space). \qed

Knowing that ${\cal H}^2(\Xi)$ is a Hilbert space, it is 
straightforward from the definition of ${\cal H}^2(\Xi)$
that the left regular representation $L$ of $G$ on  
${\cal H}^2(\Xi)$ is unitary. 
In the terminology of [FT99] or [K99] this means that 
${\cal H}^2(\Xi)$ is a {\it $G$-invariant Hilbert space of 
holomorphic functions on $\Xi$}. Those Hilbert spaces
feature the following properties:

\Remark 5.4. (a)  The inclusion mapping 
${\cal H}^2(\Xi)\into {\cal O}(\Xi)$ is continuous and so 
all point evaluations
$${\rm ev}_z\: {\cal H}^2(\Xi)\to \C, \ \ f\mapsto f(z)\qquad (z\in
\Xi)$$
are continuous. In particular for every 
$z\in \Xi$ there exists an element 
$K_z\in {\cal H}^2(\Xi)$ such that $f(z)=\la f,K_z\ra$.  
\par\nin
(b) The Hardy space admits a reproducing 
kernel
$$K_\Xi\: \Xi\times \Xi\to\C, \ \ (z,w)\mapsto \la K_w,K_z\ra\ .$$
Note that $K_\Xi$ is holomorphic in the first
and antiholomorphic in the second variable. 
We call $K_\Xi$ the {\it Cauchy-Szeg\"o kernel}
of $\Xi$. 
\par\nin
(c) By construction the representation 
$(L, {\cal H}^2(\Xi))$ of $G$ is unitary. This can also
be phrased by saying that $K_\Xi$ is $G$-invariant:
$$K_\Xi(gz,gw)=K_\Xi(z,w)$$
for all $g\in G$, $z,w\in \Xi$.
\qed
\msk

\subheadline{The boundary value map}

Using the standard semigroup techniques for Hardy spaces 
(cf.\ [H\'O\O91] or [N99, Ch.\ XIV]) one 
easily shows that the prescription 
$$L\: \oline \Gamma\to B({\cal H}^2(\Xi)), \ \ 
(L(\gamma)f)(z)\:=f(\gamma^{-1}z)$$
defines a strongly continuous involutive contractive 
representation of $\oline\Gamma$ whose restriction to 
$\Gamma$ is a holomorphic mapping. Here {\it strongly continuous} means
that the mapping $\oline \Gamma\to {\cal H}^2(\Xi),\  \gamma\mapsto
L(\gamma)f$ is continuous for all $f\in  {\cal H}^2(\Xi)$; 
{\it involutive} means $L(\gamma^*)=L(\gamma)^*$ for all 
$\gamma\in \oline\Gamma$, and {\it contractive} means that 
$\|L(\gamma)\|\leq 1$ for all $\gamma\in \oline \Gamma$. 
Note that $L$ is involutive implies in particular that 
$L\res_G$ is unitary. 

\par In order to define the boundary value map we 
first have to exhibit a good $G$-invariant 
subspace of ${\cal H}^2(\Xi)$ whose elements extend continuously 
to the boundary. A good choice herefore is the space 
of analytic vectors ${\cal H}^2(\Xi)^\omega$ of the 
representation $(L, {\cal H}^2(\Xi))$ of $G$.  
According to  [KN\'O97, App.] we have
$${\cal H}^2(\Xi)^\omega=\bigcup_{\gamma\in \Gamma}
L(\gamma){\cal H}^2(\Xi)\ .$$
In particular it follows from the compression property (5.1)
that all functions in ${\cal H}^2(\Xi)^\omega$ holomorphically 
extend over the closure $\oline \Xi$ of $\Xi$. 
As a consequence we obtain 
a well defined boundary value mapping 
$$b^\omega\: {\cal H}^2(\Xi)^\omega\to L^2(G/H), 
\ \ f\mapsto \big(gH\mapsto \lim_{\gamma\to \1
\atop \gamma\in \Gamma} f(\gamma^{-1} gz_1)\big)$$
which is $G$-equivariant and continuous.

\Theorem 5.5.  The boundary value mapping 
$$b^\omega\: {\cal H}^2(\Xi)^\omega\to L^2(G/H), 
\ \ f\mapsto \big(gH\mapsto \lim_{\gamma\to \1
\atop \gamma\in \Gamma} f(\gamma^{-1} gz_1)\big)$$
extends to a $G$-equivariant isometric embedding 
$$b\: {\cal H}^2(\Xi)\to L^2(G/H)\ .$$

\Proof. This follows from our discussion above 
and the standard semigroup techniques (cf.\ [H\'OO91] or 
[N99, Ch.\ XIV]).\qed

\subheadline{The isomorphism between ${\cal H}^2(\Xi)$
and ${\cal H}^2({\cal D})$}

\par In the sequel we will identify $\Xi$ with ${\cal D}$ under 
the biholomorphism $\Phi\: \Xi\to {\cal D}$ from Section 2. 
We let $G$ act on ${\cal O}(\Xi)$
via the left regular representation $L$ and  we 
let $S$ (or its double covering if necessary) act on ${\cal O}({\cal D})$
via the cocycle action $\pi_h$

$$(\pi_h(s)f)(z)= J_h(s^{-1},z)^{-1} f(s^{-1}z)$$
for all $s\in S$, $z\in {\cal D}$ and $f\in {\cal O}({\cal D})$. 
Then we have an isomorphism of topological vector spaces 
$$\Psi\: {\cal O}({\cal D})\to
 {\cal O}(\Xi), 
\ \ F \mapsto  \psi F$$
with inverse mapping 
$$\Psi^{-1}\: {\cal O}(\Xi)\to {\cal O}({\cal D}), 
\ \ f\mapsto {1\over \psi} f $$
Note that $\Psi$ and $\Psi^{-1}$ are defined  since $\psi$ is a
zero free holomorphic function according to Proposition 3.5.

\Lemma 5.6. The maps $\Psi$ and $\Psi^{-1}$ are 
$G$-equivariant. 

\Proof. Since $J_h=J_{-\rho_n}$ this follows from
Proposition 3.5.\qed

{}From Lemma 5.6 we obtain immediately that the mapping 

$$\Psi\: (\pi_h\res_G,{\cal H}^2({\cal D}))\to (L, {\cal O}(\Xi)), 
\ \ f\mapsto \psi f$$ 
is a continuous $G$-equivariant embedding. More precisely 
we have: 

\Theorem 5.7. The mapping 
$$\Psi\: (\pi_h\res_G,{\cal H}^2({\cal D}))\to (L, {\cal H}^2(\Xi)), 
\ \ f\mapsto \psi f$$ 
is a $G$-equivariant isomorphism of Hilbert spaces. 

\Proof.  Recall the element 
$X^0\in \z(\uu)\cap \q_*$. For $s\in ]-\infty, 0[$ define 
$$\gamma_s\:=\exp(isX^0)$$
and note that $\gamma_s\in Z(U_\C)\cap \Gamma$. Furthermore 
we have $\gamma_s^{-1} (z_1)=z_{e^s}$ which follows by simple 
$\su(1,1)$-reduction using the results in Appendix B. 
We first show that ${1\over \psi}f\in {\cal H}^2({\cal D})$ for all 
$f\in {\cal H}^2(\Xi)$ and that $\|f\|=\|{1\over \psi}f\|$. 
For every $s<0$ set 
$f_s(gH)\:=f(\gamma_s^{-1} gz_1)$ for $g\in G$. 
Now it follows from the definition of ${\cal H}^2(\Xi)$ and  
Lemma 3.11 that 

$$\eqalign{\|f\|^2 &=\lim_{\gamma\to \1\atop \gamma\in \Gamma}
\int_{G/H} |f(\gamma^{-1} gz_1)|^2\ d\mu_{G/H}(gH)\cr 
&=\sup_{s<0} \int_{G/H} |f(\gamma_s^{-1} gz_1)|^2\ d\mu_{G/H}(gH)\cr 
&= 
\sup_{s<0} \int_{G/H} |f_s(gz_1)|^2\ d\mu_{G/H}(gH)\cr 
&= 
\sup_{s<0} \int_{\partial_s{\cal D}} |f_s(z)|^2|\psi(z)|^{-2}
\ d\mu_{\partial_s{\cal D}}(z)\cr 
&= 
\sup_{0\leq r<1} \int_{\partial_s{\cal D}}\left|{f(rz)\over \psi(z)}\right|^2\ 
d\mu_{\partial_s{\cal D}}(z). \cr}\leqno(5.4)$$ 
Now consider the function 

$$F\: [0,1]\times \partial_s{\cal D}, \ \ (r,z)\mapsto 
\left|{\psi(z)\over\psi(rz)}\right|\ .$$  
Note that $F$ is continuous, hence bounded since defined on a compact 
set. In particular we find a constant $C>0$ such that 
${1\over |\psi(z)|}\geq C {1\over |\psi(rz)|}$ for all $r,z$. 
With  this information 
we obtain from (5.4) that 

$$\|f\|^2 \geq C \sup_{0\leq r<1} 
\int_{\partial_s{\cal D}}\left|{f(rz)\over \psi(rz)}\right|^2\ 
d\mu_{\partial_s{\cal D}}(z),\leqno(5.5) $$ 
i.e., ${1\over \psi} f\in {\cal H}^2({\cal D})$. 
Moreover from $F(1,z)=1$ for all $z\in \partial_s{\cal D}$
we now obtain from (5.4) and (5.5) that 
$$\eqalign{\|f\|^2 &=\lim_{r\to 1}
\int_{\partial_s{\cal D}}\left|{f(rz)\over \psi(z)}\right|^2\ 
d\mu_{\partial_s{\cal D}}(z)\cr 
&=\lim_{r\to 1}\int_{\partial_s{\cal D}}|F(r,z)|^2\,
\left|{f(rz)\over \psi(rz)}\right|^2\ 
d\mu_{\partial_s{\cal D}}(z).\cr }$$
Thus we have shown that $\Psi^{-1}$ maps 
${\cal H}^2(\Xi)$ isometrically into ${\cal H}^2({\cal D})$. 
Finally, reversing the arguments 
from above, also yields that $\Psi$ maps ${\cal H}^2({\cal D})$ 
isometrically into ${\cal H}^2(\Xi)$, concluding the proof of the 
theorem.\qed

An immediate conclusion from Theorem 5.7 in 
conjunction with Theorem 5.5 and Theorem 4.2 now is:

\Theorem 5.8. The image of the boundary value map $b: {\cal H}^2(\Xi)\to 
L^2(G/H)$ is a mulitiplicity free subspace 
of $L^2(G/H)_{\rm mc}$  of full spectrum.\qed

\subheadline{The Cauchy-Szeg\"o kernel}

We conclude this section with the derivation 
of  the Cauchy-Szeg\"o kernel 
for ${\cal H}^2(\Xi)$.

\Theorem 5.9. If $K_h$ is the reproducing kernel
 of ${\cal H}^2({\cal D})$
then the Cauchy-Szeg\"o kernel $K_\Xi$ is given by 
$$K_\Xi(z,w)= \psi (z)\oline {\psi (w)} K_h(z,w). $$

\Proof.  In this proof denote by $(\cdot| \cdot)$ the scalar 
product on ${\cal H}^2(\Xi)$ and by $\la\cdot, \cdot\ra$
the scalar product on ${\cal H}^2({\cal D})$. Fix $z\in\Xi$
and let $f\in {\cal H}^2(\Xi)$. Consider the function 
$w\mapsto \oline{\psi(z)}\psi(w) K_h(w,z)$
and note that this function lies in ${\cal H}^2(\Xi)$. 
We now have 
$$\eqalign{(f|\oline{\psi(z)}\psi(\cdot)
K_h(\cdot,z))
&=\psi(z) (f|\psi(\cdot) K_h(\cdot,z))\cr 
&=\psi(z) \la {1\over \psi(\cdot)}f,K_h(\cdot,z)\ra\cr 
&=\psi(z) {1\over \psi(z)} f(z)=f(z).}$$
This concludes the proof of the theorem.\qed

\Example 5.10. Let us consider the case of $S=G\times G$. 
We keep the notation of Example 3.6 and identify 
$\Xi={\cal D}$ with ${\cal D}_G\times {\cal D}_G^{\rm opp}$. 
Write $K_h^G(z,w)$ for the kernel of the classical 
Hardy space on ${\cal D}_G$. Then the kernel $K_h$ of the 
classical Hardy space on ${\cal D}$ is given by 

$$K_h\big((z_1,w_1), (z_2, w_2)\big)=K_h^G(z_1,z_2)\oline{K_h^G(w_1,w_2)}
$$
for all $(z_1,w_1)$, $(z_2,w_2)$ in ${\cal D}$. From Theorem 5.8 
and (3.4) we thus obtain that 
$$K_\Xi\big((z_1,w_1), (z_2, w_2)\big)={K_h^G(z_1,z_2)\oline{K_h^G(w_1,w_2)}
\over K_h^G(z_1,w_1)\oline {K_h^G(z_2,w_2)}}\ .\leqno(5.6)$$
In particular for $G=\SU(1,1)$ where ${\cal D}_G=D$ is 
the unit disc we have $K_h(z,w)={1\over 1-z\oline w}$ and so 
$$K_\Xi\big((z_1,w_1), (z_2, w_2)\big)={(1-z_1\oline{w_1})
(1-\oline{z_2} w_2)\over (1-z_1\oline{z_2})(1-\oline{w_1}w_2)}.\qeddis

\subheadline{ Final remarks}

It is possible to extend the results in this
paper to the cases in Theorem 1.12 (ii), (iii) and (iv). 
Let us explain this for groups $G$ which are 
structure groups of Euclidean Jordan algebras (the case 
in Theorem 1.12(ii)). Here one has 
$$\partial_d\Xi_0=\underbrace{\coprod_{j=1}^n G/H_j}_{\partial_s\Xi}
\amalg G/K\amalg G/K$$
with each $G/H_j$ non-compactly causal (for more 
precise information see [GK02b, Sect.\ 4]). 
For every boundary component $G/H_j$ of $\partial_d\Xi_0$ one 
can define a Hardy space 
$${\cal H}_j^2(\Xi_0)=\{ f\in {\cal O}(\Xi_0)\: 
\sup_{\gamma\in \Gamma} \int_{G/H_j} |f(\gamma^{-1} gz_j)|^2
\ d\mu_{G/H_j}(gH_j)<\infty\}\ .$$
We expect a  boundary value map
$$b_j\: {\cal H}_j^2(\Xi_0)\into L^2(G/H_j)$$
with image a full  subspace in $L^2(G/H_j)_{\rm mc}$. 
We also expect for each boundary component 
$G/H_j$ a $G$-invariant isometric isomorphism 
$${\cal H}_j^2(\Xi_0)\into {\cal H}^2({\cal D})$$
with the classical Hardy space.

\sectionheadline{ Appendix A: parameter calculations}

In this appendix we will derive some technical results on parameters
of unitary highest weight representations which are needed in Section 2. 

\par Unless otherwise specified $\s$ denots a simple Hermitian 
Lie algebra with compact Cartan subalgebra $\t\subeq \uu$. 
Write $\gamma_1, \ldots, \gamma_n$
for a maximal system of strongly orthogonal roots. 
Define 
$$\zeta\:={1\over 2}(\gamma_1+\ldots+\gamma_r)$$
and note that $\zeta\in i\z(\uu)^*$. 
Note that we can choose the positive system $\Delta^+$ such that 
$\beta=\gamma_1$ is the highest root. We write 
$\check \gamma_j\in i\t$ for the coroot of $\gamma_j$. Denote by 
$\e$ the subspace of $\t$ which is spanned 
by the $i\check \gamma_j$. Write 
$\Sigma\:=\Delta\res_\e\bs\{0\}$ for the restricted root system. 
Define $\Sigma_n\:=\Delta_n\res_\e$,
$\Sigma_c\:=\Delta_c\res_\e\bs\{0\}$.
Similarly define $\Sigma^+$, $\Sigma_n^+$, $\Sigma_c^+$. 
By results of Harish-Chandra and Moore we have: 

$$\Sigma_n^+=\{ {1\over 2}(\gamma_i+\gamma_j)\: 1\leq i, j\leq r\}$$
and 
$$\Sigma_c^+=\{ {1\over 2}(\gamma_i-\gamma_j)\: 1\leq i< j\leq r\}
\cup\{{1\over 2}\gamma_i\: 1\leq i\leq r\}\ .$$
The second term in $\Sigma_c^+$ only appears if 
$\s$ is not of tube type. 

\par Define the  number $d\:=\dim_\C \s_\C^{{1\over
2}(\gamma_i-\gamma_j)}
=\dim_\C \s_\C^{{1\over 2}(\gamma_i+\gamma_j)}$ 
for $i\neq j$. This number does not depend on $i$ and $j$, only 
on $\s$. The table is as follows: 
\gsk 
\centerline{\bf Table II}
$$\vbox{\tabskip=0pt\offinterlineskip
\def\tablerule{\noalign{\hrule}}
\halign{\strut#&\vrule#\tabskip=1em plus2em&
\hfil#\hfil&\vrule#&\hfil#\hfil &\vrule#\tabskip=0pt\cr\tablerule
&&\omit\hidewidth $\s$\hidewidth&& 
 $d$ & \cr\tablerule
&& $\sp(n,\R)$ && $1$ & \cr\tablerule
&& $\su(p,q)$ && $2$ & \cr\tablerule
&& $\so(2,n)$  && $n-2$ \ ($n>2$) & \cr\tablerule 
&& $\so^*(2n)$  && $4$ & \cr\tablerule
&& $\e_{6(-14)}$ && $8$ & \cr\tablerule
&& $\e_{7(-25)}$ && $8$ & \cr\tablerule
}}$$

\Theorem A.1. {\rm (Wallach) (cf.\ [W79])}
Let $\lambda_z\:=z\zeta$ for 
$z\in \R$. Then the set of scalar unitary highest weights 
$\lambda_z\in i\z(\uu)^*$ is parametrized by 
$$z\in {\Bbb W}\quad\hbox{ with}\quad {\Bbb W}= ]-\infty,  -{d(r-1)\over
2}[ \cup\{ -{d(r-1)\over 2}, 
-{d(r-2)\over 2}, \ldots, 0\}.\qeddis

The set ${\Bbb  W}$ in Theorem A.1 is called the {\it Wallach set}. Note 
that $\lambda_z$ is regular means precisely 
that $z< -{d(r-1)\over 2}$. 

\Lemma A.2. If $\s$ is of tube type, i.e, if $\Sigma$ is 
of type $C_r$, then the following assertions hold:

\item{(i)} $\rho_n=(1+{d(r-1)\over 2})\zeta.$
\item{(ii)} $\la \rho, \beta\ra={1\over 2}(1+d(r-1)).$
\item{(iii)} If $\lambda\in i\z(\uu)^*$ is such that 
$\la 2\lambda +\rho,\beta\ra <0$, then $\lambda$ is regular. 
\item{(iv)} The Hardy space weight $\lambda_h=-\rho_n$ satisfies 
the condition $\la 2\lambda_h+\rho, \beta\ra<0$. 

\Proof. (i) [\'O\O99, Lemma 2.5]. 
\par\nin (ii) Note that $\rho=\rho_n+\rho_c$. Thus (ii) is 
immediate from (i) and the structure of $\Sigma_c^+$. 
\par\nin (iii) Write $\lambda=\lambda_z$. Then by 
(ii) the condition $\la 2\lambda +\rho,\beta\ra<0$ means that 
$$z<-{1\over 2}(1+d(r-1))$$
and Wallach's Theorem A.1 implies the assertion. 
\par\nin (iv) We obtain from (i) that 

$$\la 2\lambda_h+\rho, \beta\ra=\la -\rho_n +\rho_c,\gamma_1\ra 
=-{1\over 2}(1+{d(r-1)\over 2}) +\la \rho_c,\gamma_1\ra\ .$$
Now $\la\rho_c,\gamma_1\ra ={d(r-1)\over 4}$. Thus 
$$\la 2\lambda_h +\rho, \beta\ra =-{1\over 2}(1+{d(r-1)\over 2})
+ {d(r-1)\over 4}=-{1\over 2}<0\ ,$$
as was to be shown. \qed 

\Proposition A.3.  Let $(\s,\g,\h)$ be a causally symmetric 
triple. Then the following assertions hold: 

\item{(i)} The Hardy space parameter $\lambda_h=-\rho_n$ is regular.
\item{(ii)} The Hardy space parameter $\lambda_h=-\rho_n$ 
satisfies the condition $\la \lambda_h +\rho_G^c,\beta\ra <0$. 

\Proof. If $\s=\g\oplus\g$ is the group case, then $\g$ is of tube 
type and $\rho_G^c ={1\over 2}\rho$. Thus (i) and (ii) follow in 
this case from Lemma A.2. Hence we may assume that $\s$ is simple. 
Since $\s$ is of 
tube type, (i) follows from Lemma A.2(i) and  Theorem A.1. 
\par It remains to check that (ii) is satisfied for the pairs
$$(\s,\g)=(\sp(2n,\R),\,  \sp(n,\R)),\,   (\su(2n,2n),\,
\sp(n,\C)),\quad {\rm and}\quad (\so(2,n), \so(1,n))\, .$$ 
\par Let us start with $(\sp(2n,\R), \sp(n,\R))$. In this 
case it is easy to see that $\la \rho_G^c,\beta\ra =
{1\over 2}\la \rho,\beta\ra$. Thus we obtain 
$$\la \lambda_h+\rho_G^c,\beta\ra ={1\over 2}\la
2\lambda_h+\rho,\beta\ra<0$$
by Lemma A.2 (iv).
\par Next consider the case $(\su(2n,2n), \sp(n,\C))$. Note 
that the restricted root system of $\sp(n,\C)$ is of type 
$C_n$ with with root space dimension $2$ for all roots. 
Thus we obtain from Lemma A.2(ii) that 
$$\la \rho_G^c,\beta\ra ={1\over 2} \big({1\over 2}(1+ 2(n-1))\big)
=-{1\over 4} +{n\over 2}\ .$$
On the other hand we have 
$$\la \lambda_h,\beta\ra =-{1\over 2}(1+ 1\cdot {2n-1\over 2})
=-({1\over 4} +{n\over 2})\ .$$
Thus $\la \lambda_h +\rho_G^c,\beta\ra =-{1\over 2}<0$. 
\par Finally consider the pair $(\so(2,n), \so(1,n))$. 
Here one has 

$$\la \lambda_h+\rho_g^c,\beta\ra =\la \lambda_h,\beta\ra +\la 
 \rho_g^c,\beta\ra =-{n\over 4} +{n-1\over 4}=-{1\over 4}<0\ ,$$
concluding the proof of the proposition. \qed

\sectionheadline{Appendix B: Structure theory for totally
symmetric triples}
 
In this appendix  we will develop some structure theory 
for causally symmetric triples $(\s,\g,\h)$. 
The results 
will be useful for certain $\su(1,1)$-reductions used 
throughout the paper. Also we will give 
a structural description of the maximal 
parabolic subgroup $P$ of $S$. 
We will state all results without proof since they entirely 
consist of minor variations of results already in the literature
(cf.\ [H\'O96, Ch.\ 1, Ch.\ 4, App.\ A]).  

\par  As in Section 2 let $\t\subeq \uu$ denote a compact $\tau$-stable  
Cartan subalgebra of $\s$. Set $\e\:=\t\cap\q_*$ and write 
$\Delta_\e=\Delta(\s_\C, \e_\C)$ for the restricted root 
system with respect to $\e_\C$. We can choose $\Delta^+$ (without
changing $\Delta_n^+$) such that 
$\Delta_\e^+\:=\Delta^+\res_\e\bs \{0\}$ is a positive system 
for $\Delta_\e$. Further we set $\Delta_{\e,n}^+\:=\Delta_n^+\res_\e$. 
Note that $\Delta_\e$ is a root 
system of type $C_r$. Write $\gamma_1,\ldots, \gamma_r$
for a maximal set of long positive  strongly orthogonal roots in 
$\Delta_\e$. Define 
$H^j\in i\e$ by $\gamma_k(H^j)=2\delta_{kj}$.
Note that 

$$X^0= -i{1\over 2}(H^1+\ldots+ H^r)\ .$$

We can choose 
$E^{\pm j}\in \s_\C^{\pm\gamma_j}$ such that 

$$\oline {E^j}= E^{-j}, \quad \tau(E^j)=E^{-j},\quad \theta(E^j)=-E^{-j}
\quad \hbox{and}\quad  [E^j, E^{-j}]=H^j.\leqno(B.1)$$
Note that the elements $\{H^j, E^j, E^{-j}\}$ form an 
$\su(1,1)$-triple in $\s_\C$.  Set 
$$Z^j\:=-i(E^j-E^{-j}) \qquad \hbox{and}\qquad Y^j:=E^j +E^{-j}$$
and note that $Z^j\in \q_*\cap\p_*$ and $Y^j\in \g\cap\p_*=\p$. 

If we set $X^j\:=-iH^j$ then the elements $\{X^j, Y^j, Z^j\}$ 
form an $\su(1,1)$-triple in 
$\s$ with 
relations

$$[X^j, Y^j]=2Z^j, \quad [X^j, Z^j]=2Y^j, 
\quad [Y^j, Z^j]= -2X^j.\leqno(B.2)$$

Our choice of the maximal abelian subspace $\a\subeq \p$ then is

$$\a\:=\bigoplus_{j=1}^r \R Y^j.\leqno(B.3)$$ 

Further 

$$\b\:=\bigoplus_{j=1}^r \R Z^j\leqno(B.4)$$
defines a maximal abelian subspace in $\q_*\cap\p_*$. 
Define elements 

$$Y^0\:={1\over 2}(Y^1+\ldots+ Y^r)\quad\hbox{and}\quad Z^0\:={1\over 2}
(Z^1+\ldots +Z^r)\ .$$
Note that $\Spec({\ad Y^0})=\Spec ({\ad Z^0})=\{-1, 0,1\}$. 
Next we define a Cayley transform 

$${\bf C}\: \s_\C\to \s_\C,\ \ X\mapsto  e^{i{\pi\over 2}\ad Y^0}(X)\ .$$
Notice the relations

$${\bf C}(H^j)=Z^j, \quad {\bf C}(Z^j)=-H^j, \quad {\bf C}(Y^j)=Y^j\ .$$

Since $\Spec({\ad Y^0})=\{-1, 0,1\}$, the prescription 
$$\sigma(X)=(e^{i\pi \ad  Y^0}\circ \theta)(X)\qquad (X\in\s)$$
defines an involution on $\s$. 
Note that $\sigma$ commutes with $\tau$ by construction. 
Write $\s=\h^\star \oplus\q^\star$ for the $\sigma$-eigenspace 
decomposition of $\s$ and $\g=\h\oplus\q$ for the 
$\sigma$-eigenspace decomposition of $\g$, i.e., $\h=\h^\star\cap\g$
and $\q=\g\cap \q^\star$. 
Then 
$$\a\subeq \p\cap\q\qquad \hbox{and}\qquad \b\subeq \q_*\cap\q^{\star}
\cap\p_*.\leqno(B.5)$$

\ssk  As in Section 2 we let $S_\C$ be a simply connected Lie group 
with Lie algebra $\s_\C$ and $S$, $G$, $U$, $K$ the analytic 
subgroups of $S_\C$ corresponding to $\s$, $\g$, $\uu$  and 
$\k$. Our final  goal is to describe the maximal parabolic 
subgroup $P={\bf C}(U_\C P^-)\cap S$ of $S$.  
For that write $\Sigma_\b=\Sigma(\s,\b)$ for the restricted root
system of $\s$ with respect to $\b$. Note that 
$\Sigma_\b={\bf C}^t(\Delta_\e)$ and so $\Sigma_\b$ is of type 
$C_r$. Define $\Sigma_{\b,n}^+\:=C^t(\Delta_{\e,n}^+)$ and write 

$$\n_\b^{\pm}\:=\bigoplus_{\alpha\in \pm \Sigma_{\b,n}^+} \s^\alpha.
\leqno(B.6)$$
Note that $\n_\b^{\pm}$ are both abelian and real forms of 
${\bf C}(\p^{\pm})$. Write $N_{\b}^\pm$ for 
the analytic subgroups of $S$ corresponding to 
$\n_\b^\pm$. 
Write $L\:={\bf C}(U_\C )\cap S$. Then $L$ is the Levi 
factor of $P$ and we have 

$$P=L\rtimes N_\b^-.\leqno(B.7)$$
We will show that $L=H^\star\:=S^\sigma$. 
First we describe $\l\:={\rm Lie}(\l)$. 
For that define subalgebras

$$\s(0)\:=\z_\s(Y^0)\quad\hbox{and}\quad \g(0)=\z_\g(Y^0)\ .$$
Then $\s(0)=\s^{\sigma\theta}$ and $\g(0)=\g^{\sigma\theta}$
and therefore 

$$\g(0)=(\h\cap\k)\oplus(\p\cap\q) \quad\hbox{and}\quad 
\s(0)=(\h^\star\cap \uu)\oplus(\q^\star\cap\p_*).\leqno(B.8)$$

Thus we obtain 

$$\eqalign{\l &={\bf C}(\uu_\C)\cap \s=(\uu_\C\cap \s(0))\oplus
\big({\bf C}(\uu_\C)\cap (\h^\star\cap \p_*\oplus
\q^\star\cap\uu)\big)\cr 
&=(\h^\star \cap \uu)\oplus (\h^\star\cap \p_*)=\h^\star.\cr}$$
Hence $L=H^\star$.  
Summarizing our discussion we have shown: 

\Lemma B.1. Let $P={\bf C}(U_\C P^-)\cap S$ be the stabilizer 
of $x_1\in \partial_s{\cal D}$ in $S$. Then the Levi-decomposition 
of the maximal parabolic subgroup $P$ is given by 

$$P=H^\star\rtimes N_\b^-.\qeddis

\def\entries{

\[AG90 Akhiezer, D.\ N., and S.\ G.\ Gindikin, {\it On Stein
extensions of 
real symmetric spaces}, 
Math.\ Ann.\ {\bf 286}, 1--12, 1990

\[BS97 van den Ban, E., and H.\ Schlichtkrull, {\it The most
continuous part of the Plan\-che\-rel decomposition for a reductive 
symmetric space}, Ann. of Math. {\bf (2) 145} (1997), no. {\bf 2},
267--364

\[BS01a  --- , {\it The Plancherel decomposition for a
reductive symmetric space I. The spherical part},  preprint 2001 

\[BS01b ---,  {\it The Plancherel decomposition for a reductive symmetric
space II. Representation theory}, preprint  2001

\[Ber96 Bertram, W., {\it On some causal and conformal groups}, J. Lie 
Theory {\bf 6} (1996), 215--247

\[Bet97 Betten, F., {\it Causal compactification of compactly causal 
spaces}, TAMS, to appear 

\[B\'O01 Betten, F., and G. \'Olafsson, {\it Causal compactification
and Hardy spaces for spaces of Hermitian type}, Pacific, J. of
Math. Vol. {\bf 200 (2)} (2001),  273--312

\[D{\'O}Z01 Davidson, M., G. \'Olafsson, G. Zhang, {\it Segal-Bargmann transform
on Hermitian symmetric spaces and orthogonal polynomials}. Preprint, 2001

\[De98  Delorme, P., {\it Formule de Plancherel pour les espaces 
sym\'etriques r\'eductifs},  Ann. of Math. {\bf (2) 147}
(1998), no. {\bf 2}, 417--452

\[EHW83 Enright, T., R. Howe, and N.\ Wallach, {\it A classification of 
unitary highest weight modules},  Representation theory of reductive 
groups (Park City, Utah, 1982), 97--143, Progr. Math., {\bf 40}, 
Birkh\"auser Boston, Boston, MA, 1983

\[FT99 Faraut, J., and  E. G. F. Thomas, {\it Invariant Hilbert spaces
of holomorphic functions}, J. Lie Theory {\bf 9} (1999), no. {\bf 2}, 
383--402

\[GG77 Gel'fand, M., and S.G. Gindikin, {\it Complex Manifolds 
whose Skeletons are real Lie Groups, and Analytic Discrete Series 
of Representations}, Funct. 
Anal. and Appl. {\bf 11}(1977), 19--27

\[G98 Gindikin, S., {\it Tube domains in Stein symmetric spaces}, 
Positivity in Lie theory: open problems, 81--97, de Gruyter Exp. Math., 
{\bf 26}, de Gruyter, Berlin, 1998

\[GK02a Gindikin, S., and B.\ Kr\"otz, {\it Complex crowns of Riemannian 
symmetric spaces and non-compactly causal symmetric spaces}, 
Trans. Amer. Math. Soc. {\bf 354} (2002), no. {\bf 8}, 
3299--3327

\[GK02b ---, {\it Invariant Stein domains in Stein symmetric spaces
and a non-linear complex convexity theorem}, IMRN {\bf 18} (2002), 959--971

\[GM01 Gindikin, S., and T. Matsuki, {\it Stein Extensions
of Riemann Symmetric Spaces and Dualities of Orbits on Flag 
Manifolds}, Transformation groups, to appear 

\[GrK01 Groenevelt, W., and E.\ Koelnik, {\it Meixner functions
and polynomials related to Lie algebra representations}, preprint

\[H\'O96 Hilgert, J.\ and 
G.\ \'Olafsson, ``Causal Symmetric Spaces, Geometry and
Harmonic Analysis,'' Acad. Press, 1996 

\[H\'O\O91  Hilgert, J., G.\ 'Olafsson, and B.\ \O rsted, {\it 
Hardy spaces on affine symmetric spaces}, 
J. Reine Angew. Math. {\bf 415} (1991), 189--218

\[K99  Kr\"otz, B., {\it The Plancherel theorem for biinvariant 
Hilbert spaces},  Publ. Res. Inst. Math. Sci. {\bf 35} (1999), 
no. {\bf 1}, 91--122

\[K01 ---, {\it Formal dimension for semisimple symmetric spaces}, 
Compositio Math. {\bf 125} (2001), no. {\bf 2}, 155--191

\[KN02 Kr\"otz, B., and K.--H. Neeb, {\it Unitary spherical 
highest weight representations}, Trans. Amer. Math. Soc. {\bf 354} (2002), no. {\bf 3}, 
1233--1264

\[KN\'O97 Kr\"otz, B., K.--H. Neeb, and G.\ \'Olafsson, {\it 
Spherical representations and mixed symmetric spaces}, 
Representation Theory {\bf 1} (1997), 424--461

\[K\'O02  Kr\"otz, B., and G.\ \'Olafsson, {\it The c-function for a
non-compactly causal symmetric space}, Invent. math. {\bf 149(3)} (2002), 647--659

\[KS01a Kr\"otz, B., and R.J. Stanton, {\it Holomorphic extension of
representations: (I)  
automorphic functions}, Annals of Mathematics, to appear

\[KS01b Kr\"otz, B., and R.J. Stanton, {\it Holomorphic extension of
representations: (II)  
geometry and harmonic analysis}, preprint

\[N99 Neeb, K.--H., ``Holomorphy and Convexity in Lie Theory,'' 
Expositions in Mathematics, de Gruyter {\bf 28}, 1999

\[Ne99 Neretin, Y. A., {\it Matrix analogs of the} $B${\it -function
and Plancherel formula for Berezin Kernel representation}.
Math. Sb. {\bf 191} (2000), 57--100. Translation
in Sb. Math. {\bf 191} (2000), 683--715

\[\'O00  \'Olafsson, G., {\it  Analytic continuation in representation 
theory and harmonic analysis},  Glo\-bal analysis and harmonic analysis 
(Marseille-Luminy, 1999), 201--233, S\'emin. Congr., {\bf 4}, 
Soc. Math. France, Paris, 2000

\[\'O\O88 \'Olafsson, G., and B.\ \O rsted, {\it  The holomorphic
discrete series for affine symmetric 
spaces. I}, J. Funct. Anal. {\bf 81} (1988), no. {\bf 1}, 126--159

\[\'O\O96 ---, {\it Generalizations of the Bargmann transform},
Lie theory and its applications in physics (Clausthal, 1995)
(H.-D. Doebner, V.K. Dobrev, and J. Hilgert eds.)
World Scientific, River Edge, New Jersey, 1996

\[\'O\O99 ---, {\it  Causal compactification 
and Hardy spaces},  Trans. Amer. Math. Soc. {\bf 351} (1999), no. {\bf
9}, 
3771--3792

\[O82 Olshanski, G. I., {\it Invariant cones in Lie algebras, 
Lie semigroups, and the holomorphic discrete series}, Funct. Anal. and Appl. 
{\bf 15}, 275--285 (1982)

\[S86 Stanton, R. J., {\it Analytic Extension of the holomorphic discrete 
series}, Amer. J. Math. {\bf 108} (1986), 1411--1424

\[W79 Wallach, N., {\it The analytic continuation of the 
discrete series.I, II}, Trans. Amer. Math. Soc. {\bf 251} (1979), 
1--17, 19--37

\[Z01 Zhang, G., {\it Berezin transform on real bounded symmetric domains},
Trans. Amer. Math. Soc., {\bf 353} (2001), 3769--3787

}

{\sectionheadline{\bf References}
\frenchspacing
\entries\par}
\dlastpage 
\bye
\end